\documentclass[final,leqno]{siamltex}
\usepackage{graphicx}
\usepackage{algorithm} 
\usepackage{algorithmic} 

\title{On Interpolation Approximation: Convergence rates  for polynomial interpolation for functions of limited regularity\thanks{This work was
supported by National Science Foundation of China (No. 11371376).}}

\author{Shuhuang Xiang\thanks{School of
Mathematics and Statistics, Central South University, Changsha, Hunan
410083, P. R. China.}}

\begin{document}
\maketitle

\begin{abstract}
The convergence rates on polynomial interpolation in most cases are estimated by Lebesgue constants. These estimates may be overestimated for some special points of sets for functions of limited regularities. In this paper, by applying the Peano kernel theorem and  Wainerman's lemma, new formulas on the convergence rates are considered. Based upon these new estimates, it shows that the interpolation at strongly normal pointsystems can achieve the optimal convergence rate, the same as the best polynomial approximation. Furthermore, by using the asymptotics on Jacobi polynomials,  the convergence rates are established for Gauss-Jacobi, Jacobi-Gauss-Lobatto or Jacobi-Gauss-Radau pointsystems. From these results, we see that the interpolations at the Gauss-Legendre, Legendre-Gauss-Lobatto pointsystem, or at strongly normal pointsystems, has essentially the same approximation accuracy compared with those at the two Chebyshev piontsystems, which also illustrates the equally accuracy of the Gauss and Clenshaw-Curtis quadrature. In addition,  numerical examples illustrate the perfect coincidence with the estimates, which means the convergence rates are optimal.
\end{abstract}



\pagestyle{myheadings}
\thispagestyle{plain}

\section{Introduction}
%

A central problem in approximation theory
is the construction of simple functions that are easily implemented on computers and approximate well a given set of functions.

There exist many investigations for the behavior of continuous functions approximated by polynomials. Weierstrass \cite{Weierstrass} in 1885 proved the well known result that every continuous function $f(x)$ in $[-1,1]$ can be uniformly approximated as closely as desired by a polynomial function. This result has both practical and theoretical relevance, especially in polynomial interpolation.

Polynomial interpolation is a fundamental tool in many areas of scientific computing. Lagrange interpolation is a classical technique for approximation of continuous functions. Let us denote by
\begin{equation}\label{fundamentalpoints}
-1\le x_n^{(n)}<x_{n-1}^{(n)}<\cdots<x_2^{(n)}<x_1^{(n)}\le 1
\end{equation}
the $n$ distinct points in the interval $[-1,1]$ and let $f(x)$ be a function
defined in the same interval. The $n$th Lagrange interpolation
polynomial of $f(x)$ is unique  and given by the formula
\begin{equation}\label{lagrange}
L_n[f] =\sum_{k=1}^nf(x_k^{(n)})\ell_k^{(n)}(x),\quad \ell_k^{(n)}(x)=\frac{\omega_n(x)}{\omega'_n(x_k^{(n)})(x-x_k^{(n)})},
\end{equation}
where $\omega_n(x)=(x-x_1^{(n)})(x-x_2^{(n)})\cdots(x-x_n^{(n)})$.

There is a well developed theory that quantifies the convergence or divergence of the Lagrange interpolation polynomials (Brutman \cite{Brutman1978,Brutman1997} and Trefethen \cite{Trefethen2}).  Two key notions for interpolation in a given set of points are that of the \textit{Lebesgue function}
\begin{equation}\label{lebesguefun}
   \lambda_n(x)=\sum_{k=1}^n\big|\ell_k^{(n)}(x)\big|
\end{equation}
and \textit{Lebesgue constant}
\begin{equation}\label{lebesguecon}
  \Lambda_n=\max_{x\in[-1,1]} \lambda_n(x),
\end{equation}
which are of fundamental importance (Cheney \cite{Cheney}, Davis \cite{Davis} and Szeg\"{o} \cite{Szego}).
The Lebesgue constant can also be interpreted as the $\infty$-norm of the projection operator $L_n: C([-1,1])\rightarrow {\cal P}_{n-1}$
 $$\Lambda_n = \sup_{f} {\|L_n[f]\|_{\infty}\over \|f\|_{\infty}}, $$ where ${\cal P}_{n-1}$  is the set of polynomials of degree less than or equal to $n-1$.

Based upon the Lebesgue constant, the interpolation error can be estimated by
\begin{equation}\label{best}
\big\|L_n[f]-f\big\|_{\infty}\le (1+\Lambda_n)\big\|p^*_{n-1}-f\big\|_{\infty},
\end{equation}
where $p^*_{n-1}$ is the best polynomial approximation of degree $n-1$. Thus, the Lebesgue constant $\Lambda_n$ indicates  how good the interpolant $L_n[f]$ is in comparison with the best polynomial approximation $p^*_{n-1}$.

The study of the Lebesgue constant $\Lambda_n$ originated more than $100$ years ago. Comprehensive reviews can be found in Brutman \cite{Brutman1997}, Lubinsky \cite{Lubinsky}, Trefethen \cite[Chapter 15]{Trefethen2}, etc. For an arbitrarily given
system of points $\{x_1^{(n)},x_2^{(n)},\ldots,x_n^{(n)}\}_{n=1}^{\infty}$, Bernstein \cite{Bernstein1914} and Faber \cite{Faber} in 1914 obtained that
$$
\Lambda_n\ge \frac{1}{12}\log n,
$$
which, together with the boundedness principle,  implies that there exists a  continuous function $f(x)$ in $[-1,1]$ for which the sequence $L_n[f]$ ($n=1,2,\ldots$) is not uniformly convergent to $f$ in $[-1,1]$\footnote{Gr\"{u}nwald \cite{Grunwald1935} in 1935 and Marcinkiewicz \cite{Marcinkiewicz} in 1937, independently, showed that even for the Chebyshev points of first kind
$$
x_k^{(n)}=\cos\left(\frac{2k-1}{2n}\pi\right),\quad k=1,2,\ldots,n,\quad n=1,2,\ldots,
$$
there is a continuous  function $f(x)$ in $[-1,1]$ for which the sequence $L_n[f]$ is divergent everywhere in $[-1,1]$.}. More precisely, Erd\"{o}s \cite{Erdos1961} and Brutman \cite{Brutman1978} proved that
\begin{equation}\label{LB}
\Lambda_n\ge \frac{2}{\pi}\log n +C \mbox{\, for some constant $C$ (\cite{Erdos1961})};\quad \Lambda_n\ge  \frac{2}{\pi}\left(\gamma_0+\log \frac{4}{\pi}\right)+\frac{2}{\pi}\log n \, (\cite{Brutman1978}),
\end{equation}
where $\gamma_0=0.577\ldots$ is the Euler's constant. In particular, for equidistant pointsystem $$\left\{x_k^{(n)}=-1+\frac{2k}{n-1}\right\}_{k=0}^{n-1},$$ Sch\"{o}nhage \cite{Schonhage1961} showed that
$$
\Lambda_n\sim \frac{2^n}{e(\log (n-1)+\gamma_0)(n-1)},\quad n\rightarrow \infty.
$$
Additionally, Trefethen and Weideman \cite{Trefethen3} established that
$$
\frac{2^{n-3}}{(n-1)^2}\le \Lambda_n\le \frac{2^{n+2}}{n-1},\quad n\ge 0.
$$
Then generally, the set of equally spaced points is a bad choice for Lagrange interpolation (see Runge \cite{Runge1901}).

Whereas, for well chosen sets of points, the growth of $\Lambda_n$ may be extremely slow as $n\to\infty$:
\begin{itemize}
\item Chebyshev pointsystem of first kind $T_n=\left\{x_k^{(n)}=\cos\left(\frac{2k-1}{2n}\pi\right)\right\}_{k=1}^{n}$: An asymptotic estimate of $\Lambda_n(T_n)$ was given by Bernstein \cite{Bernstein} as
\begin{equation}
\Lambda_n(T_n)\sim \frac{2}{\pi}\log n, \quad n\rightarrow \infty,
\end{equation}
which is improved by Ehlich and Zeller \cite{Ehlich1966}, Rivlin \cite{Rivlin1974} and Brutman \cite{Brutman1978} as
$$
\frac{2}{\pi}\left(\gamma_0+\log \frac{4}{\pi}\right)+\frac{2}{\pi}\log n<\Lambda_n(T_n)\le 1+ \frac{2}{\pi}\log n,\quad n=1,2,\ldots.
$$

\item Chebyshev pointsystem of second kind $U_n=\left\{x_k^{(n)}=\cos\left(\frac{k}{n-1}\pi\right)\right\}_{k=0}^{n-1}$ (also called Chebyshev extreme  or Clenshaw-Curtis points \cite{Trefethen1}): Ehlich and Zeller \cite{Ehlich1966} proved that
\begin{equation}\label{ChBY2}
\Lambda_n(U_n)=\left\{\begin{array}{ll}\Lambda_{n-1}(T_{n-1}),&n=2,4,6,\ldots\\
\Lambda_{n-1}(T_{n-1})-\alpha_n,\quad 0\le \alpha_n<\frac{1}{(n-1)^2},&n=3,5,7,\ldots.\end{array}\right.
\end{equation}

\item The roots of Jacobi polynomial $P_n^{(\alpha,\beta)}(x)$ ($\alpha,\beta>-1$):  The asymptotic estimate of $\Lambda_n(J_n)$ was found by Szeg\"{o} \cite{Szego} as
\begin{equation}\label{Jacobi}
  \Lambda_n(J_n)=\left\{\begin{array}{ll}O(n^{\gamma+\frac{1}{2}}),&\gamma>-\frac{1}{2}\\
  O(\log n),&\gamma\le-\frac{1}{2}\end{array},\right.\quad \gamma=\max\{\alpha,\beta\}.
\end{equation}

\end{itemize}
Comparing Equations (1.7), (1.8) and (1.9) with (1.6), we see that the two Chebyshev pointsystems and the Jacobi pointsystem with $\gamma\le -\frac{1}{2}$ are nearly optimal and of order $O(\log n)$.

Nevertheless, it is worth noting that if $f(x)$ has an absolutely continuous $(k-1)$st
derivative $f^{(k-1)}$ on $[-1,1]$ for some $k\ge 1$ and its $k$-th derivative $f^{(k)}$ is of  bounded variation ${\rm Var}(f^{(k)})<\infty$, Mastroianni and Szabados \cite{Mastroianni1995}, Trefethen \cite{Trefethen2} and Xiang et al. \cite{XiangChenWang} proved that
\begin{equation}\label{apprerr}
\|f-L_n[f]\|_{\infty}=O(n^{-k}),
\end{equation}
where $L_n[f]$ is at the $n$ Chebyshev points of first or second kind, which has the same asymptotic order as
$\|f-p^*_{n-1}\|_{\infty}$ for the best approximation $p_{n-1}^*$, following de la Vall\'{e}e Poussin \cite{Poussin1908}. In particular, for $f(x)=|x|$, the error on the $L_n[f]$  at the above two Chebyshev pointsystems satisfies
$$
\|f-L_n[f]\|_{\infty}\le \frac{4}{\pi(n-1)}
$$
(see \cite{Trefethen2,XiangChenWang}), while
$$\|f-p^*_{n-1}\|_{\infty}\sim \frac{\beta}{n}, \quad 0.2801685<\beta<0.2801734
$$
(see Bernstein \cite{Bernstein1914b} and Varga and Capenter \cite{Varga1985}). Thus,  the error estimate (1.5) by using the Lebesgue constant may be overestimated for some special points of sets for functions of limited regularities.

Moreover, it has been observed, by Clenshaw-Curtis \cite{Clenshaw1960} and O'Hara and Smith \cite{Hara1968}, that $n$-point Gauss quadrature and $n$-point Clenshaw-Curtis quadrature have essentially the same accuracy, which has been showed  recently by Trefethen \cite{Trefethen1,Trefethen2}, Brass and Petras \cite{BrassPetras} and Xiang and Bornemann \cite{XiangBornemann}. Both of these two quadrature
are derived from the interpolation polynomial $L_n[f]$ by
 $$
Q_n[f]=\int_{-1}^1L_n[f](x)dx,
$$
 based on the $n$ Gauss-Legender and Clenshaw-Curtis points, respectively. From this observation, we may conclude that the corresponding interpolation $L_n[f]$ based on  these two pointsystems may have the same convergence rate. However, it can not be derived from (1.5).

In this paper, we present new convergence rates of the interpolation polynomials for functions of limited regularities, based upon the famous Peano kernel theorem \cite{Peano1913} and applying an interesting Wainerman's lemma \cite{Wainerman}. Suppose $f(x)$ has an absolutely continuous $(r-1)$st
derivative $f^{(r-1)}$ on $[-1,1]$, and  its $r$-th derivative $f^{(r)}$ is of  bounded variation ${\rm Var}(f^{(r)})<\infty$. We will show that
\begin{equation}\label{gerror}
\|f-L_n[f]\|_{\infty}\le \frac{\pi^r {\rm Var}(f^{(r)})}{(n-1)(n-2)\cdots(n-r)}\max_{1\le j\le n}\|\ell_j^{(n)}\|_{\infty},
\end{equation}
which leads to
\begin{equation}
\|f-L_n[f]\|_{\infty}=O(n^{-r}\max_{1\le j\le n}\|\ell_j^{(n)}\|_{\infty}).
\end{equation}
The Lebesgue constant $\Lambda_n=\max_{x\in[-1,1]}\sum_{k=1}^n\big|\ell_k^{(n)}(x)\big|$ is replaced by  $\max_{1\le j\le n}\|\ell_j^{(n)}\|_{\infty}$ in some sense since $\|f-p^*_{n-1}\|_{\infty}=O(n^{-r})$ \cite{Poussin1908}.

Particularly, from (1.12), it directly follows that the interpolation $L_n[f]$ at a strongly normal pointsystem (see Fej\'{e}r \cite{Fejer1916}) can achieve the optimal convergence rate as $O(\|f-p^*_{n-1}\|_{\infty})$.

Furthermore, $\|\ell_j\|_{\infty}$ can be explicitly estimated for Gauss-Jacobi, Jacobi-Gauss-Lobatto or Jacobi-Gauss-Radau pointsystems, by using the asymptotics on Jacobi polynomials given by Szeg\"{o} \cite{Szego} and some results given in Kelzon \cite{Kelzon,Kelzon2}, V\'{e}rtesi \cite{Vertesi1980,Vertesi1983b}, Sun \cite{Sun}, Prestin \cite{Prestin}, Kvernadze \cite{Kvernadze}, Vecchia et al. \cite{Vecchia}, etc., as follows
\begin{itemize}
\item For the $n$ Gauss-Jacobi points:
$$
\max_{1\le j\le n}\|\ell_j^{(n)}\|_{\infty}=O(n^{\max\{\gamma-\frac{1}{2},0\}}),\quad \gamma=\max\{\alpha,\beta\}.
$$

\item For the $n$ Jacobi-Gauss-Lobatto points (the roots of $(1-x^2)P_{n-2}^{(\alpha,\beta)}(x)=0$):
$$
\max_{1\le j\le n}\|\ell_j^{(n)}\|_{\infty}=\left\{\begin{array}{ll}
O\left(n^{-\min\{0,\alpha+\frac{1}{2},\beta+\frac{1}{2}\}}\right),&-1<\alpha,\beta\le \frac{3}{2}\\
O\left(n^{-\min\{0,\alpha+\frac{1}{2},2+\alpha-\beta,\frac{5}{2}-\beta\}}\right),&-1<\alpha\le \frac{3}{2},\beta> \frac{3}{2}\\
O\left(n^{-\min\{0,\beta+\frac{1}{2},2+\beta-\alpha,\frac{5}{2}-\alpha\}}\right),&\alpha>\frac{3}{2}, -1<\beta\le\frac{3}{2}\\
O\left(n^{-\min\{0,2+\alpha-\beta,2+\beta-\alpha,\frac{5}{2}-\alpha,\frac{5}{2}-\beta\}}\right),&\alpha,\beta>\frac{3}{2}\end{array}.\right.
$$

\item For the $n$ Jacobi-Gauss-Radau points $(1-x)P_{n-1}^{(\alpha,\beta)}(x)$
$$
\max_{0\le j\le n-1}\|\ell_j^{(n)}\|_{\infty}=\left\{\begin{array}{ll}O\left(n^{-\min\{0,\alpha+\frac{1}{2},\alpha-\beta\}}\right),&-1<\alpha\le\frac{1}{2}\\
O\left(n^{-\min\{0,\frac{1}{2}-\beta,\frac{5}{2}-\alpha,\alpha-\beta\}}\right),&\alpha>\frac{1}{2}\end{array}.\right.
$$

\item For the $n$ Jacobi-Gauss-Radau points $(1+x)P_{n-1}^{(\alpha,\beta)}(x)$
$$
\max_{1\le j\le n}\|\ell_j^{(n)}\|_{\infty}=\left\{\begin{array}{ll}O\left(n^{-\min\{0,\beta+\frac{1}{2},\beta-\alpha\}}\right),&-1<\beta\le\frac{1}{2}\\
O\left(n^{-\min\{0,\frac{1}{2}-\alpha,\frac{5}{2}-\beta,\beta-\alpha\}}\right),&\beta>\frac{1}{2}\end{array}.\right.
$$
\end{itemize}
From the above estimates, we see that the interpolation at the Gauss-Legendre or at the Legendre-Gauss-Lobatto pointsystem, has essentially the same approximation accuracy compared with those at the two Chebyshev piontsystems. All of them satisfy that $\max_{1\le j\le n}\|\ell_j^{(n)}\|_{\infty}=O(1)$ (for more general cases see {\sc Fig}. 1.1). In addition, the convergence rate is attainable illustrated by some functions of limited regularities.

\begin{figure}[htbp]
\centerline{\includegraphics[height=4.5cm,width=4.5cm]{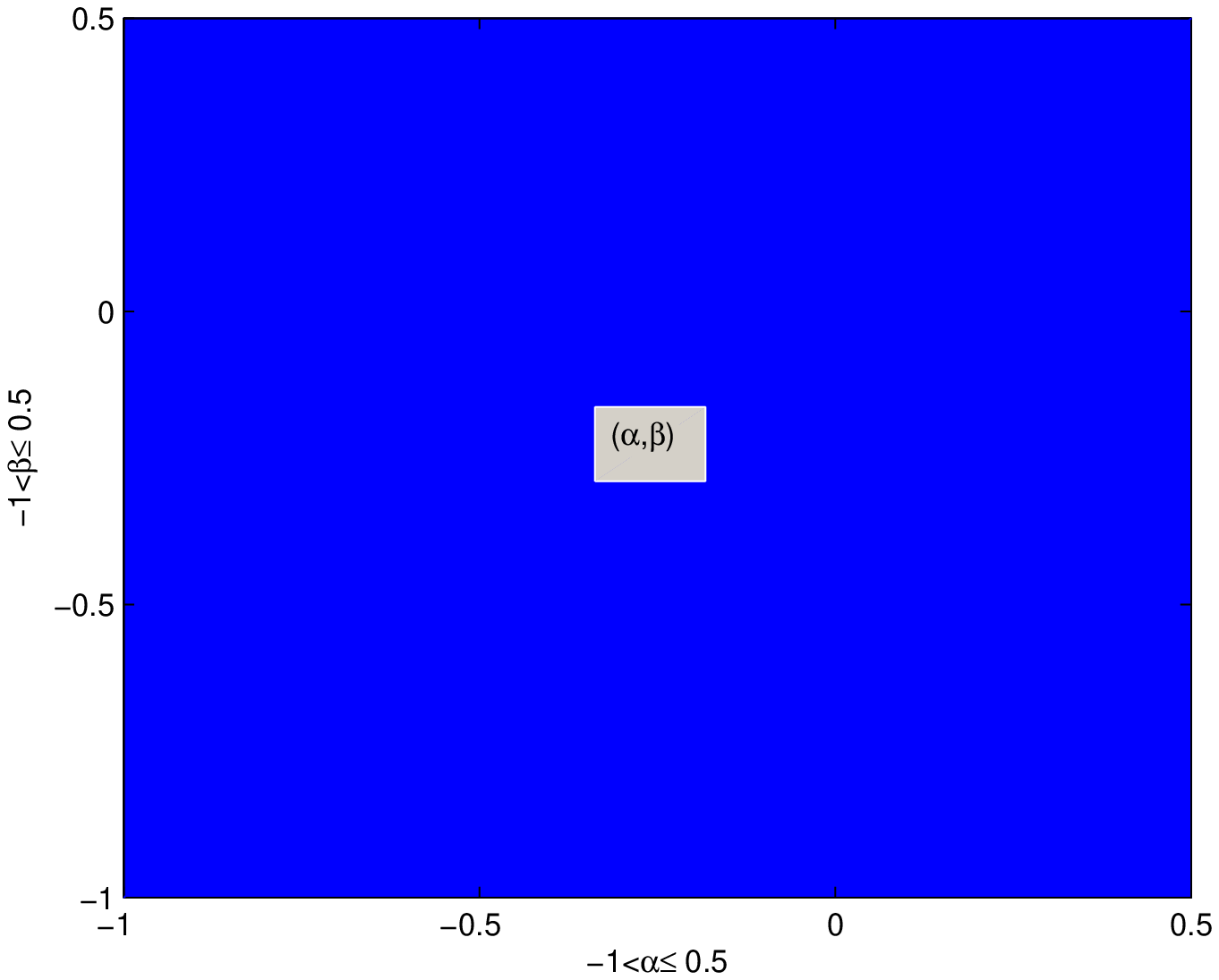}\hspace{0.6cm}\includegraphics[height=4.5cm,width=4.5cm]{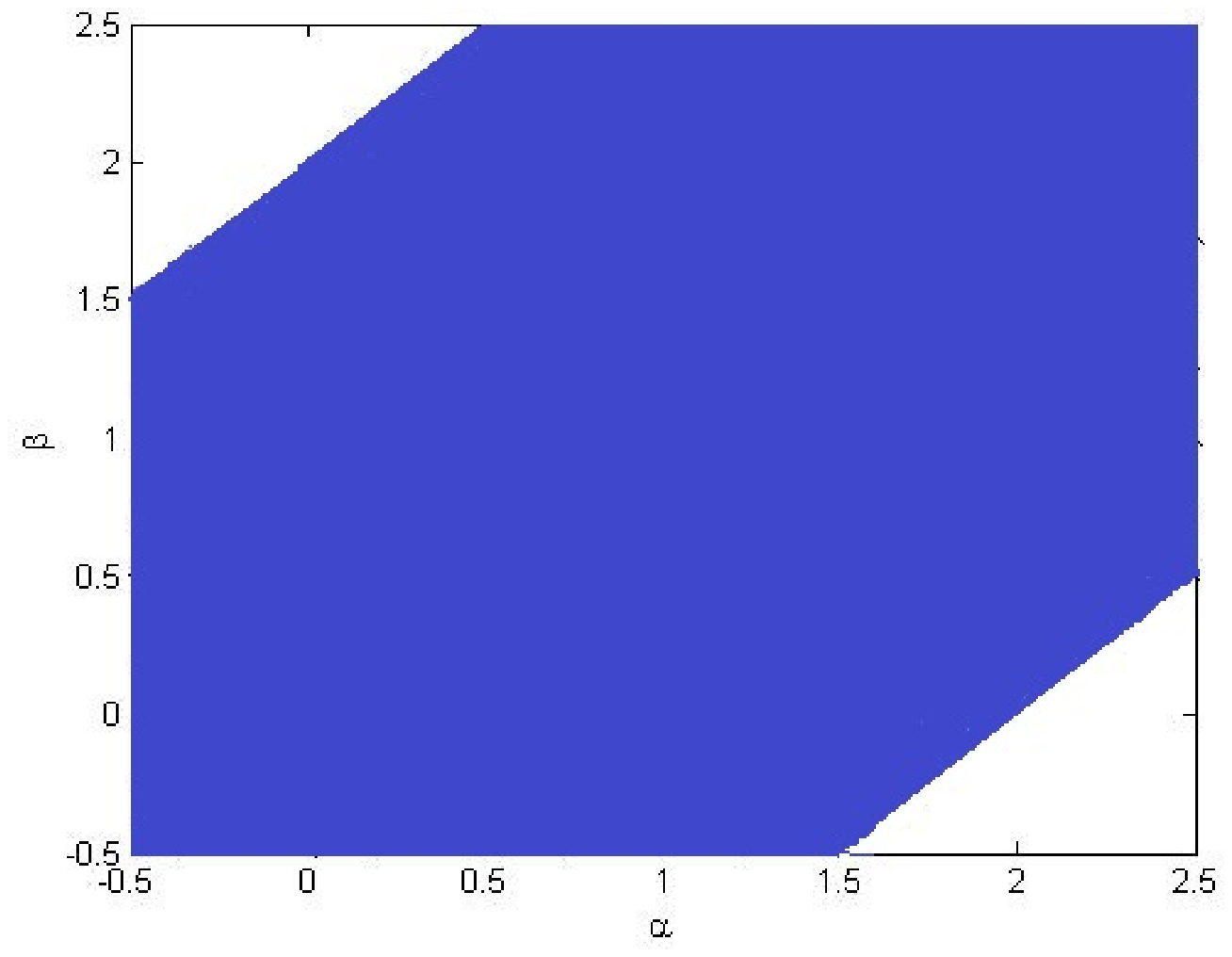}\hspace{0.6cm}
\vspace{1cm}\includegraphics[height=4.5cm,width=4.5cm]{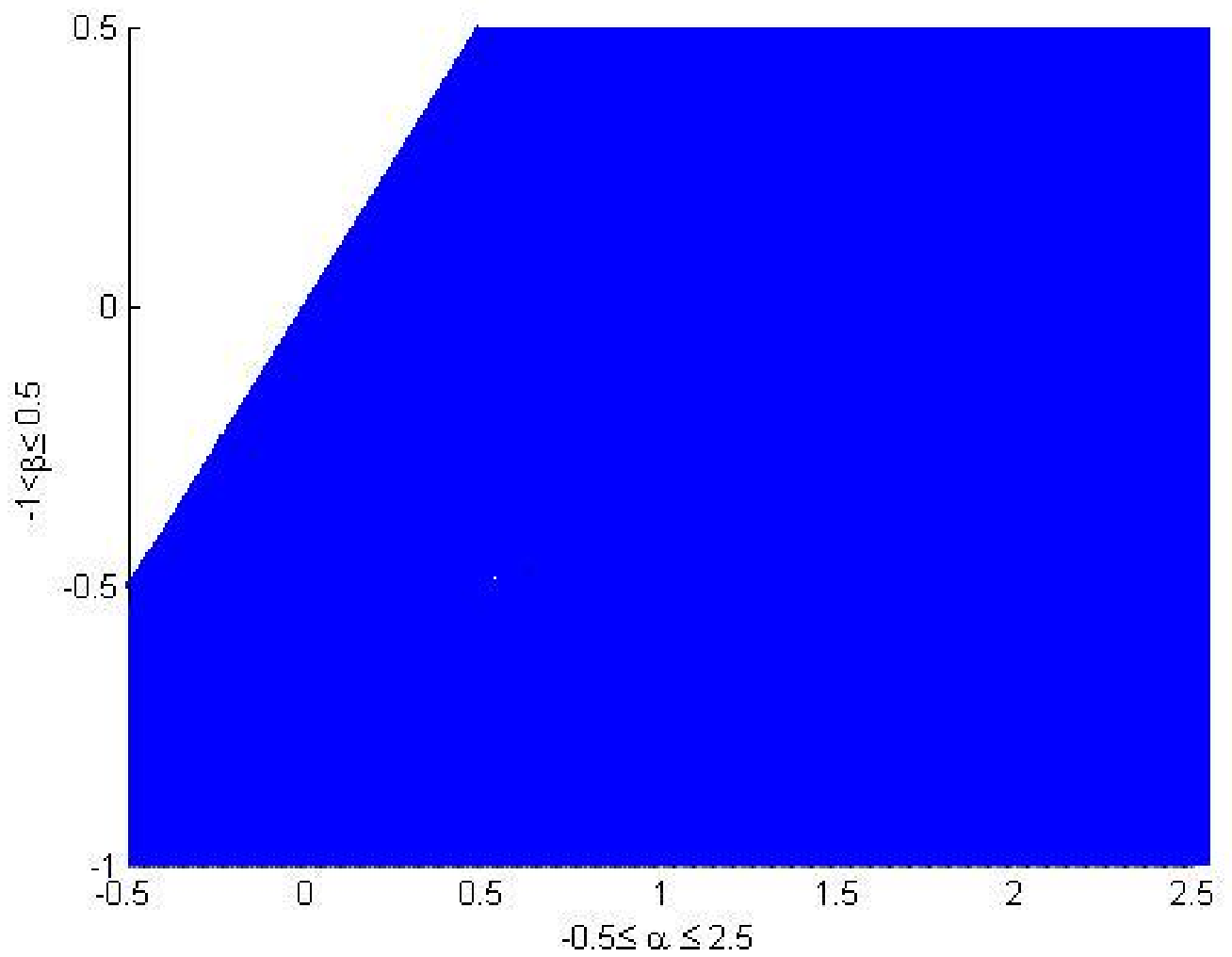}}
\caption{The neighbourhood on $(\alpha,\beta)$ such that $\|f-L_n[f]\|_{\infty}=O(n^{-r})$ for the  Gauss-Jacobi pointsystems (left), Jacobi-Gauss-Lobatto pointsystems (middle) and Jacobi-Gauss-Radau ($(1-x)P_{n-1}^{(\alpha,\beta)}(x)=0$) pointsystems (right), respectively: {\rm Var}$(f^{(r)})<\infty$.}
\end{figure}

Thus, the best approximation polynomial is challenged by the interpolation polynomials at the special pointsystems showed in {\sc Fig.} 1.1. Furthermore, we will see that the interpolation polynomials at the special pointsystems perform much better than the best approximation polynomial for approximation the derivatives $f'$ and $f''$ by $L_n'[f]$, $L_n''[f]$,  $[p_{n-1}^*]'$ and  $[p_{n-1}^*]''$, respectively, illustrated by numerical examples in the final section.

It is worthy of special mention that the interpolation polynomial $L_n[f]$, at the  Gauss-Jacobi, Jacobi-Gauss-Lobatto or  Gauss-Jacobi-Radau pointsystem,
 can be efficiently evaluated by applying the second barycentric  formula
 $${\displaystyle
L_n[f](x)=\frac{\sum_{j=1}^n\frac{\lambda_j}{x-x_j}f(x_j)}{\sum_{j=1}^n\frac{\lambda_j}{x-x_j}}},
$$
which is robust in the presence of rounding
errors \cite{Higham2004} and costs overall computational complexity  $O(n)$  \cite{Berrut2004}, where the nodes $x_j$ and the barycentric weights $\lambda_j$ are computed by \textbf{jacpts} and the formulas given in \cite{Hale2012,WangXiang2012,Wang2013}, respectively.  A {\sc Matlab} routine \textbf{jacpts}, which uses the algorithm in \cite{HaleTownsend} for the computation of these nodes and weights, can be found in {\sc Chebfun} system \cite{Chebfun}. For more details on this topic, see Salzer \cite{Salzer}, Henrich \cite{Henricibook1982}, Berrut and Trefethen \cite{Berrut2004}, Higham \cite{Highambook,Higham2004}, Glaser et al. \cite{Glaser}, Wang and Xiang \cite{WangXiang2012}, Bogaert et al. \cite{Bogaert}, Hale and Trefethen \cite{Hale2012}, Hale and Townsend \cite{HaleTownsend}, Trefethen \cite{Trefethen2}, Wang et. al. \cite{Wang2013}, etc. {\sc Matlab} routines can be found in {\sc Chebfun}
system \cite{Chebfun} and Xiang and He \cite{XiangHe}.

The paper is organized as follows: In section 2, we present the error of $f(x)-L_n[f](x)$ for each fixed $x\in [-1,1]$ by using the Peano representation and the bounded variation. In section 3, we introduce the interesting Wainerman's lemma and deduce the error bound on $\|f-L_n[f]\|_{\infty}$ by $\max_{1\le j\le n}\|\ell_j^{(n)}\|_{\infty}$.
We consider, in section 4, the estimates of $\|\ell_j^{(n)}\|_{\infty}$ and derive the convergence rates for the interpolation polynomial at strongly normal pointsystems, Gauss-Jacobi, Jacobi-Gauss-Lobatto  and Jacobi-Gauss-Radau pointsystems, respectively, where the convergence rates  and attainability are illustrated by numerical experiments.


Throughout this paper,
$A\sim B$ means that there exist positive constants $C_1$ and $C_2$ such that
$$
C_1 B\le A\le C_2 B.
$$
For simplicity, in the following we abbreviate $x_k^{(n)}$ as $x_k$ and $\ell_k^{(n)}(x)$ as $\ell_k(x)$.

All the numerical results in this paper are carried out by using {\sc Matlab} R2012a on a  desktop (2.8 GB RAM, 2 Core2 (32 bit)
processors at 2.80 GHz) with Windows XP operating system.

\section{The Peano kernel theorem}
There are two general methods for deriving strict error bounds (Dahlquist and Bj\"{o}rck \cite{Dahlquist}). One applies the norms and distance formula together with the Lebesgue constants, which
often overestimates the error. The other is due to the Peano kernel theorem. 

 Suppose ${\cal L}$ a continuously linear functional that maps functions $f\in C([-1,1])$ to $R$ satisfying ${\cal L}(f_1+f_2) = {\cal L} f_1 + {\cal L}f_2$ for any $f_1, f_2 \in C([-1,1])$ and ${\cal L}(\alpha f) = \alpha {\cal L}f$ for any scalar $\alpha$. In addition, we assume ${\cal L}[{\cal P}_{r-1}] = \{0\}$
for some $r\in \{1, 2,\ldots\}$, where ${\cal P}_{r-1}$ denotes  the set of polynomials with degree less than or equal to  $r-1$.

The Peano kernel theorem (Peano \cite{Peano1913}, see also Kowalewski \cite{Kowalewski1932}, Schmidt
\cite{Schmidt1935} and Mises \cite{Mises1936}) is the identity
\begin{equation}\label{Peano}
{\cal L}[f]=\int_{-1}^1f^{(r)}(t)K_r(t)dt
\end{equation}
holding for all such functions $f\in C^{r}([-1,1])$, where $K_r(t)=\frac{1}{(r-1)!}{\cal L}[(x-t)_+^{r-1}]$ and
\begin{equation}\label{frac}
(x-t)_+^{r-1}=\left\{\begin{array}{ll}
(x-t)^{r-1},&x\ge t\\
0,&x<t\end{array}\right.\quad(r\ge 2),\quad (x-t)_+^0=\left\{\begin{array}{ll}
1,&x\ge t\\
0,&x<t.\end{array}\right.\quad(r=1).
\end{equation}

For each fixed $x\in [-1,1]$, we consider the special functional ${\cal L} = E_n$, where $E_n[f](x)$ is defined for $\forall f\in C([-1,1])$ by $$E_n[f](x)=f(x)-\sum_{j=1}^nf(x_j)\ell_j(x)=f(x)-L_n[f](x)$$ with $-1\le x_n<x_{n-1}<\cdots<x_2<x_1\le 1$. $E_n[f]$ is a continuously linear functional
 since $|E_n[f](x)-E_n[g](x)|\le (1+\Lambda_n)\|f-g\|_{\infty}$ for arbitrary $f,g\in C([-1,1])$, and then  by the Peano theorem \cite{Peano1913} $E_n[f]$ can be represented if $f\in C^{r}([-1,1])$ for $n\ge r$ as
\begin{equation}\label{PeanoThm}
E_n[f](x)=\int_{-1}^1f^{(r)}(t)K_r(t)dt
\end{equation}
with
\begin{equation}
K_r(t)=\frac{1}{(r-1)!}(x-t)_+^{r-1}-\frac{1}{(r-1)!}\sum_{j=1}^n(x_j-t)_+^{r-1}\ell_j(x).
\end{equation}
Particularly, from (2.3) it implies
$$
|E_n[f](x)| \le \|f^{(r)}\|_{\infty}\int_{-1}^1|K_r(t)|dt\le 2\|f^{(r)}\|_{\infty}\|K_r\|_{\infty}.
$$

Similar to the Peano kernel for quadrature \cite{BrassPetras}, the kernel for interpolation satisfies the following proposition.

\begin{proposition}(Peano representation)
Let
\begin{equation}
K_s(t)=\frac{1}{(s-1)!}(x-t)_+^{s-1}-\frac{1}{(s-1)!}\sum_{j=1}^n(x_j-t)_+^{s-1}\ell_j(x),\quad s=1,2,\ldots.
\end{equation}
Then for $s\ge 2$, the Peano kernel satisfies $K_{s}(-1)=K_{s}(1)=0$ and can be rewritten as
\begin{equation}\label{rep}
K_s(u) =\int_u^1K_{s-1}(t)dt,\quad s=2,3,\ldots.
\end{equation}
\end{proposition}
\begin{proof}
From the definition of $K_s$ in (2.5), it is easy to verify that $K_s(-1)=K_s(1)=0$ by using $\sum_{j=1}^{n}\ell_j(t)\equiv 1$ for $t\in [-1,1]$. Furthermore,
we find that
\begin{equation}\,\,\,\,
\begin{array}{lll}
{\displaystyle\frac{1}{(s-2)!}\int_u^1(x-t)_+^{s-1}dt}&=&\left\{\begin{array}{ll}
0,& u>x\\
{\displaystyle\frac{1}{(s-2)!}\int_u^x(x-t)^{s-1}dt=\frac{1}{(s-1)!}(x-u)^{s-1}},&u\le x\end{array}\right.\\
&=&{\displaystyle\frac{1}{(s-1)!}(x-u)_+^{s-1}.}\end{array}
\end{equation}
Define $x_0=1$ and $x_{n+1}=-1$ and suppose $x_{m+1}< u\le x_m$ for some nonnegative integer $m$. By (2.7), similarly  we have
\begin{equation}\quad \begin{array}{lll}
{\displaystyle \frac{1}{(s-2)!}\sum_{j=1}^n\int_u^1(x_j-t)_+^{s-2}\ell_j(x)dt}&=&{\displaystyle\frac{1}{(s-2)!}\sum_{j=1}^m\ell_j(x)\left\{\begin{array}{ll}
0,& u>x_j\\
\int_u^{x_j}(x_j-t)^{s-2}dt,&u\le x_j\end{array}\right.}\\
&=&{\displaystyle\frac{1}{(s-1)!}\sum_{j=1}^n(x_j-u)_+^{s-1}\ell_j(x).}\end{array}
\end{equation}
Then from
$$
\int_u^1K_{s-1}(t)dt={\displaystyle\frac{1}{(s-2)!}\int_u^1(x-t)_+^{s-2}dt-\frac{1}{(s-2)!}\sum_{j=1}^n\int_u^1(x_j-t)_+^{s-2}\ell_j(x)dt},
$$
we get $\int_u^1K_{s-1}(t)dt=K_s(u)$ by (2.7) and (2.8).
\end{proof}

\vspace{0.6cm}
In the following, we consider  functions of limited regularities as
\begin{equation}\quad\quad
\begin{array}{l}\mbox{\emph{Suppose that $f(t)$ has an absolutely continuous $(r-1)$st
derivative $f^{(r-1)}$ on $[-1,1]$}}\\
\mbox{for some $r\ge 1$ with $f^{(r-1)}(t)=f^{(r-1)}(-1)+\int_{-1}^tg(y)dy$,
\emph{where $g$ is absolutely}}\\
\mbox{\emph{integrable and of bounded variation ${\rm Var}(g)<\infty$  on $[-1,1]$.}}
\end{array}
\end{equation}

From Stein and Shakarchi \cite[p. 130]{Stein2} and Tao \cite[pp. 143-145]{Tao}, we see that a function $G: [-1,1]\rightarrow R$ is absolutely continuous if and only if it takes the form $G(t)=\int_{-1}^tg(y)dy+C$ for some absolutely integrable $g: [-1,1]\rightarrow R$ and a constant $C$. It is obvious that such $g$ is not unique. Then in this paper,  we suppose $f(t)$ satisfies (2.9) and define
$$
V_r=\inf\left\{{\rm Var}(g)\,\,\Big |\,\,\begin{array}{l} f^{(r-1)}(t)=f^{(r-1)}(-1)+\int_{-1}^tg(y)dy \mbox{\, for all $t\in [-1,1]$ with $g$ being}\\
\mbox{ absolutely integrable and of bounded variation}\end{array} \right\}.
$$

{\sc Remark 1}. Here, we use the condition ``$f^{(r-1)}(t)=f^{(r-1)}(-1)+\int_{-1}^tg(y)dy$,
where $g$ is absolutely integrable and of bounded variation ${\rm Var}(g)<\infty$'' instead of ``$f^{(r)}$ is of bounded variation $V_r={\rm Var}(f^{(r)})<\infty$'' in  \cite{Trefethen1,Trefethen2}. If $f^{(r)}$  is of bounded variation, then $f^{(r+1)}$ exists almost everywhere and $f^{(r+1)}\in L^1([-1,1])$ (see Lang \cite{Lang1997} and Rudin \cite{Rudin}). Whereas, $f^{(r)}$ in  \cite{Trefethen1,Trefethen2} denotes an equivalent representation in the sense of almost everywhere. An example for $f(x)=|t|$ is given in  \cite{Trefethen1,Trefethen2}, where $f(t)$ is not differentiable at $t=0$, but $f'$ can be chosen as $$f'(t)=\left\{\begin{array}{ll}1,& t> 0\\c,&t=0\\ -1,& t< 0 \end{array}\right.,$$ then ${\rm Var}(f')=\left\{\begin{array}{ll}2,& |c|\le 1\\|1+c|+|1-c|,&{\rm otherwise}\end{array}\right.$. Using the new condition, we see that $|t|$ can be represented as  $|t|=1+\int_{-1}^tg(y)dy$ with $g(y)=\left\{\begin{array}{ll}1,& y> 0\\c,&y=0\\ -1,& y< 0 \end{array}\right.$ and $V_1=2$ is unique.


\begin{theorem}
Suppose $f(t)$ satisfies (2.9), then for $n\ge r$, we have
\begin{equation}\label{error11}
\|E_n[f]\|_{\infty}\le V_r\|K_{r+1}\|_{\infty}.
\end{equation}
\end{theorem}
\begin{proof}
Applying the Peano theorem implies that for each fixed $x\in [-1,1]$,
$$
E_n[f](x)=\int_{-1}^1f^{(s)}(t)K_{s}(t)dt,\quad s=1,2,\ldots,r-1.
$$
Then, directly following Brass and Petras \cite{BrassPetras}, integrating by parts and using $K_{r}(-1)=K_r(1)=0$ yields
$$
E_n[f](x)=\int_{-1}^1f^{(r-1)}(t)K_{r-1}(t)dt=\int_{-1}^1g(t)K_r(t)dt.
$$
Since $g$ can be written as $g=g_1-g_2$ with $g_1$ and $g_2$ are monotonically increasing, and ${\rm Var}(g)={\rm Var}(g_1)+{\rm Var}(g_2)$ (see Lang \cite[pp. 280-281]{Lang1997}). Without loss of generality, assume $g$ is monotonically increasing. Then by the second mean value theorem of integral calculus, it follows from $K_{r+1}(-1)=\int_{-1}^1K_{r}(t)dt=0$ that there exists a $\xi\in [-1,1]$ such that
$$
E_n[f](x)=g(-1)\int_{-1}^{\xi}K_r(t)dt+g(1)\int_{\xi}^1K_r(t)dt=(g(1)-g(-1))K_{r+1}(\xi)={\rm Var}(g)K_{r+1}(\xi),
$$
which leads to the desired result.
\end{proof}

\begin{lemma}\cite[Lemma 5.7.1]{BrassPetras}
Assume that
$$
\sup_{-1\le t\le 1}w(t)\sqrt{1-t^2}<\infty,\quad t_u(y)=\left\{\begin{array}{ll}0,&y<u\\
1,& y\ge u\end{array}.\right.
$$
Then, for every positive integer $\ell$ and every $u\in [-1,1]$, there is a $q_u\in {\cal P}_{\ell}$ satisfying
$$
q_u(y)\ge t_u(y) \quad \mbox{for all $y\in [-1,1]$}
$$
and
$$
\int_{-1}^1\left[t_u(y)-q_u(y)\right]w(y)dy\ge -\frac{\pi}{\ell+1}\sup_{-1\le t\le 1}w(t)\sqrt{1-t^2}.
$$
\end{lemma}

\begin{lemma}
\begin{equation}\label{peano}
|K_{s+1}(u)|\le \frac{\pi}{n-s+1}\sup_{-1\le t\le 1}|K_s(t)|.
\end{equation}
\end{lemma}
\begin{proof}In Lemma 2.3, letting $\ell=n-s-1$, $w(t)\equiv 1$, representing $q_u$ as $q_u(t)=p_{n-1}^{(s)}(t)$, and noting that $E_n[{\cal P}_{n-1}]=0$, by Theorem 2.2 we have $$0=E_n[p_{n-1}]=\int_{-1}^1p_{n-1}^{(s)}(t)K_{s}(t)dt=\int_{-1}^1q_u(t)K_{s}(t)dt.$$ Consequently, by Lemma 2.3 we get  that
$$\begin{array}{lll}
\Big|K_{s+1}(u)\Big|=\Big|\int_u^1K_s(t)dt\Big|=\Big|\int_{-1}^1K_s(t)t_u(t)dt\Big|&=&\Big|\int_{-1}^1K_s(t)\left[t_u(t)-q_u(t)\right]dt\Big|\\
&\le&{\displaystyle \frac{\pi}{n-s}\sup_{-1\le t\le 1}|K_s(t)|.}\end{array}
$$
\end{proof}

From Theorem 2.2 and Lemma 2.4 we obtain that
\begin{theorem}
Suppose $f(t)$ satisfies (2.9), then for $n\ge r+1$
\begin{equation}\label{error1}
\|E_n[f]\|_{\infty}\le \frac{\pi^r V_r}{(n-1)(n-2)\cdots(n-r)}\|K_{1}\|_{\infty}.
\end{equation}
\end{theorem}


\section{Wainerman's lemma}
In the following, we shall focus on the estimate of $\|K_{1}\|_{\infty}$.

Notice that $\sum_{j=1}^{n}\ell_j(t)\equiv 1$ for $t\in [-1,1]$ and
\begin{equation}\label{K1}
K_1(u)=(x-u)_+^{0}-\sum_{j=1}^n(x_j-u)_+^{0}\ell_j(x).
\end{equation}
If $x_1<u\le 1$, we have $K_1(u)=1$ for $u\le x$, and $K_1(u)=0$ for $u> x$. While
for $-1<u\le x_n$, we have $K_1(u)=0$ for $u\le x$, and $K_1(u)=-1$ for $u> x$. Thus, in these cases we obtain
\begin{equation}
|K_1(u)|\le 1\le \max_{1\le j\le n}\|\ell_j\|_{\infty}
\end{equation}
 since $\ell_j(x_j)=1$ for $j=1,2,\ldots,n$.

Suppose that $x_{m+1}<u\le x_{m}$ for some positive integer $m$, then for $u\le x$ we get
\begin{equation}
K_1(u)=1-\sum_{j=1}^{n}(x_j-u)_+^{0}\ell_j(x)=1-\sum_{j=1}^{m}\ell_j(x)=\sum_{j=m+1}^n\ell_j(x),
\end{equation}
while for $u> x$ we have
\begin{equation}
K_1(u)=-\sum_{j=1}^{n}(x_j-u)_+^{0}\ell_j(x)=-\sum_{j=1}^{m}\ell_j(x).
\end{equation}

\begin{lemma} (Wainerman's lemma \cite{Wainerman})
Suppose  $x_{m+1}<u\le x_m$ for some positive integer $m$, and let
$$
a_k(u)=\left\{\begin{array}{ll}
{\displaystyle\sum_{j=1}^k\ell_j(u)},&k=1,2,\ldots,m\\
{\displaystyle\sum_{j=k}^n\ell_j(u)},&k=m+1,m+2,\ldots,n\end{array}\right.
$$
and $a_0(u)=a_{n+1}(u)\equiv 0$.
Then it follows for $x_{m+1}<u< x_m$ that
\begin{equation}\label{peanokernel}{\displaystyle
{\rm sgn}(a_k(u))={\rm sgn}(\ell_k(u))=\left\{\begin{array}{ll}
(-1)^{m-k},&k=1,2,\ldots,m\\
(-1)^{k-m-1},&k=m+1,m+2,\ldots,n\end{array}\right.}
\end{equation}
and for $x_{m+1}<u\le x_m$ that
\begin{equation}\label{peanokerne2}
|a_k(u)|\le |\ell_k(u)|,\quad k=1,2,\ldots,n,
\end{equation}
where ${\rm sgn}$ denotes the sign function.
\end{lemma}
\begin{proof}
The interesting result and its proof is published in Russian in \cite{Wainerman}. For convenience and completeness, we present the proof here.

For $x_{m+1}<u<x_m$, from the definition of $\ell_k(t)$ we see that
$$\begin{array}{lll}{\displaystyle
{\rm sgn}(\ell_{k}(u))}&=& {\displaystyle {\rm sgn}\left(\frac{(u-x_1)\cdots (u-x_{k-1})(u-x_{k+1})\cdots (u-x_n)}{(x_{k}-x_1)\cdots (x_{k}-x_{k-1})(x_{k}-x_{k+1})\cdots (x_{k}-x_n)}\right)}\\
&=&{\displaystyle (-1)^{1-k}{\rm sgn}\left((u-x_1)\cdots (u-x_{k-1})(u-x_{k+1})\cdots (u-x_n)\right)},\end{array}
$$
which directly leads to the desired result (3.5) for ${\rm sgn}(\ell_k(u))$ based on $k\le m$ or $k>m$, respectively.

In the following, we will show that  ${\rm sgn}(a_k(u))$ also satisfies (3.5).

\textbf{In the case $k\le m$}: Since
$$
a_k(x_j)={\displaystyle\sum_{i=1}^k\ell_i(x_j)}=\left\{\begin{array}{ll}
1,&j=1,2,\ldots,k\\
0,&j=k+1,k+2,\ldots,n\end{array},\right.
$$
then by the Rolle's theorem it follows
$$
a'_k(y_j)=0
$$
for some $y_j$ satisfying $x_{j+1}<y_j<x_j$ for $j=1,\ldots,k-1,k+1,\ldots,n-1$.

\vspace{-0.5cm}
\begin{center}
$$
\unitlength=0.11cm
\begin{picture}(120,36)
\put(0,16){\line(1,0){120}}

\put(0,18){$-1$} \put(120,18){$1$} \put(10,18){$x_n$}\put(20,18){$x_{n-1}$}\put(46,18){$x_{k+2}$}\put(59,18){$x_{k+1}$}\put(70,18){$x_k$}\put(82,18){$x_{k-1}$}
\put(100,18){$x_{2}$}\put(115,18){$x_{1}$}

\put(0,15){$\bullet$} \put(120,15){$\bullet$} \put(10,15){$\bullet$}\put(20,15){$\bullet$}\put(46,15){$\bullet$}\put(60,15){$\bullet$}\put(70,15){$\bullet$}\put(82,15){$\bullet$}
\put(100,15){$\bullet$}\put(115,15){$\bullet$}\put(95,15){$\cdots$}\put(30,15){$\cdots$}

\put(61,30){$a_k(x_k)=1$}\put(76,30){$a_k(x_{k-1})=1$}
\put(95,30){$a_k(x_2)=1$}\put(111,30){$a_k(x_1)=1$}

\put(3,12){$a(x_n)=0$}\put(18,12){$a(x_{n-1})=0$}\put(38,12){$a(x_{k+2})=0$}\put(57,12){$a(x_{k+1})=0$}

\put(70,27){$\ast$}\put(82,27){$\ast$}
\put(100,27){$\ast$}\put(115,27){$\ast$}

\put(15,15){$\circ$}\put(55,15){$\circ$}\put(40,15){$\circ$}\put(76,15){$\circ$}\put(90,15){$\circ$}
\put(110,15){$\circ$}
\put(13,20){$y_{n-1}$}\put(53,20){$y_{k+1}$}\put(40,20){$y_{k+2}$}\put(75,20){$y_{k-1}$}\put(90,20){$y_{k-2}$}
\put(110,20){$y_{1}$}

\end{picture}
$$
\end{center}

\vspace{-0.5cm}
Note that $a_k(t)$ is a polynomial of degree $n-1$, then $a_k'(t)$ is a polynomial of degree $n-2$, which implies that $y_j$ are the exact zeros of $a'_k(t)$ and then $a'_k(t)$ has the form of
\begin{equation}
a'_k(t)=C(t-y_1)\cdots(t-y_{k-1})(t-y_{k+1})\cdots(t-y_{n-1})
\end{equation}
for some non-zero constant $C$. In addition, from (3.7) $a'_k(t)$ has alternative sign between these roots.
Then, by $a_k(x_{k+1})=0$ and $a_k(x_{k})=1$, it yields
$$
a_k'(t)>0,\quad t\in (y_{k+1},y_{k-1})
$$
and
$$
{\rm sgn}(a_k(t))=1,\quad  t\in (x_{k+1},x_{k})\subset (y_{k+1},y_{k-1})
$$
since $a(t)$ is strictly increasing in  $(x_{k+1},x_{k})$ and $a_k(x_{k+1})=0$.

 By the alternative property of $a'_k(t)$ between these roots, it deduces that
${\rm sgn}(a_k'(t))=(-1)^{j-k}$ for $t\in (y_{j+1},y_{j})$ and $j>k$,
particularly,
$$
{\rm sgn}(a_k'(t))=1,\quad t\in (y_{k+1},x_{k+1})\subset (y_{k+1},y_{k-1})
$$
and
$$
{\rm sgn}(a_k'(t))=-1,\quad t\in (x_{k+2},y_{k+1})\subset (y_{k+2},y_{k+1}),
$$
which, together with $a_k(x_{k+1})=a_k(x_{k+2})=0$, derives  ${\rm sgn}(a_k(t))=-1$ for $t\in (x_{k+2},x_{k+1})$. Similarly, applying
$${\rm sgn}(a_k'(t))=(-1)^{j-k}, \quad t\in (y_{j+1},x_{j+1});\quad\quad {\rm sgn}(a_k'(t))=(-1)^{j-k+1},\quad t\in (x_{j+2},y_{j+1}),$$
 together with $a_k(x_{j+2})=a_k(x_{j+1})=0$ for $j>k$, derives  ${\rm sgn}(a_k(t))=(-1)^{j-k+1}$ for $j>k$ and $t\in (x_{j+2},x_{j+1})$ by induction. So we get $a_k(u)=(-1)^{m-k}$.

\textbf{In the case $k>m$}: By
$$
a_k(x_j)={\displaystyle\sum_{i=k}^n\ell_i(x_j)}=\left\{\begin{array}{ll}
0,&j=1,2,\ldots,k-1\\
1,&j=k,k+1,\ldots,n\end{array},\right.
$$
applying similar arguments derives $a_k(u)=(-1)^{k-m-1}$ for $k>m$.

Furthermore, from (3.5) and the definition of $a_k(t)$, we see that immediately: for $k\le m$ and $x_{m+1}<u<x_m$,
$$
|a_k(u)|=|\ell_k(u)+a_{k-1}(u)|=|\ell_k(u)|-|a_{k-1}(u)|\le |\ell_k(u)|,
$$
and for $k>m$
$$
|a_k(u)|=|\ell_k(u)+a_{k+1}(u)|=|\ell_k(u)|-|a_{k+1}(u)|\le |\ell_k(u)|.
$$

The special case of (3.6) for $u=x_m$ directly follows from $|a_k(x_m)|=|\ell_k(x_m)|$ by the definitions of $a_k(u)$ and $\ell_k(u)$.
\end{proof}

Theorem 2.5 together with (3.2), (3.3), (3.4) and (3.6) leads to the following estimate.

\begin{theorem}
Suppose $f(t)$ satisfies (2.9), then for $n\ge r+1$
\begin{equation}\label{error12}
\|E_n[f]\|_{\infty}\le \frac{\pi^r V_r}{(n-1)(n-2)\cdots(n-r)}\max_{1\le j\le n}\|\ell_j\|_{\infty}.
\end{equation}
\end{theorem}

In the next section, we shall focus on estimates of  $\max_{1\le j\le n}\|\ell_j\|_{\infty}$ for special points of sets.


\section{Estimates of $\|\ell_j\|_{\infty}$ and convergence rates on $\|f-L_n[f]\|_{\infty}$}
For any convergent quadrature derived from polynomial interpolation at the grid points (1.1) for
$$
\int_{-1}^1f(x)w(x)dx=\int_{-1}^1f(x)d\sigma(x)
$$
for each $\sigma(x)$ of bounded variation and any analytic function $f(x)$ on $[-1,1]$,
the
clustering of the $n$ points  has a limiting  Chebyshev distribution
$$
\mu(t)=\frac{1}{\pi}\int_{-1}^t\frac{1}{\sqrt{1-x^2}}dx
$$
(see Krylov \cite[Theorem 7, p. 263]{Krylov1962}); that is, the
clustering will be asymptotically the same: on $[-1, 1]$, $n$ points will be distributed with
density
$$\frac{n}{\pi \sqrt{1-x^2}}$$
as $n$ tends to infinity (see Hale and Trefethen \cite{Hale2008} and Trefethen \cite{Trefethen1}).

Moreover, the clustering of optimal pointsystems for polynomial interpolation implies near endpoints $\pm1$ (see Z. Ditzian and V. Totik \cite{Ditzian1987} and \cite{Trefethen1}). (The Gauss-Jacobi type pointsystems have this proposition.) The density of the zeros of orthogonal polynomials has been extensively studied in Erd\"{o}s and Tur\'{a}n \cite{Erdos1938,Erdos1940}, Gatteschi \cite{Gatteschi87} and Szeg\"{o} \cite{Szego}.

\subsection{Strongly normal pointsystems}
One of the proofs of Weierstrass¡¯ approximation theorem using interpolation
polynomials was first presented by Fej\'{e}r \cite{Fejer1916} in 1916  based on the Chebyshev pointsystem of first kind $\left\{x_k=\cos\left(\frac{2k-1}{2n}\pi\right)\right\}_{k=1}^{n}$: If $f\in C([-1,1])$, then there is a unique polynomial $H_{2n-1}(f,t) $ of degree at most
$2n - 1$ such that $\lim_{n\rightarrow \infty}\|H_{2n-1}(f)-f\|_{\infty}=0$, where $H_{2n-1}(f,t)$ is determined by
\begin{equation}
H_{2n-1}(f,x_k)=f(x_k),\quad H_{2n-1}'(f,x_k)=0,\quad k = 1, 2,\ldots, n.
\end{equation}
This polynomial is known as the Hermite-Fej\'{e}r interpolation polynomial.

The convergence result has been extended to general Hermite-Fej\'{e}r interpolation of $f(x)$  at nodes (1.1) by Gr\"{u}nwald  \cite{Grunwald1942} in 1942, upon strongly normal pointsystems introduced in Fej\'{e}r \cite{Fejer1932a}: Given, respectively,
the function values $f(x_1)$, $f(x_2)$, $\ldots$, $f(x_n)$ and derivatives
$d_1$, $d_2$, $\ldots$, $d_n$ at these grids, the general Hermite-Fej\'{e}r interpolation polynomial $ H_{2n-1}(f)$ has the form of \begin{equation}\label{hermite-fejer}
\quad \,\, H_{2n-1}(f,t) =\sum_{k=1}^nf(x_k)h_k(t)+\sum_{k=1}^nd_kb_k(t),
\end{equation}
where $h_k(t)=v_k(t)\left(\ell_k(t)\right)^2$, $b_k(t)=(t-x_k)\left(\ell_k(t)\right)^2$ and
\begin{equation}
 v_k(t)=1-(t-x_k)\frac{\omega_n''(x_k)}{\omega_n'(x_k)}\mbox{\quad (see Fej\'{e}r \cite{Fejer1932b}).}
\end{equation}

The pointsystem (1.1) is called  strongly normal if for all $n$
\begin{equation}\label{stronglynormal}
 v_k(t)\ge c>0,\quad k=1,2,\ldots,n,\quad t\in [-1,1]
\end{equation}
for some positive constant $c$. The pointsystem (1.1) is called  normal if for all $n$
\begin{equation}\label{normal}
 v_k(t)\ge 0,\quad k=1,2,\ldots,n,\quad t\in [-1,1].
\end{equation}

Fej\'{e}r \cite{Fejer1932a} (also see Szeg\"{o} \cite[p. 339]{Szego}) showed that for the zeros of Jacobi polynomial $P_n^{(\alpha,\beta)}(t)$ of degree $n$ ($\alpha>-1$, $\beta>-1$)
\begin{equation}
{\small \quad \quad v_k(t)\ge \min\{-\alpha,-\beta\}\mbox{\quad for $-1<\alpha\le 0$, $-1<\beta\le 0$, $k=1,2,\ldots,n$ and $t\in [-1,1]$}.}
\end{equation}
While for the Legendre-Gauss-Lobatto pointsystem (the roots of {\small$(1-t^2)P_{n-2}^{(1,1)}(t)=0$}),
\begin{equation}
 v_k(t)\ge 1,\quad k=1,2,\ldots,n,\quad t\in [-1,1]\,\, \mbox{(Fej\'{e}r \cite{Fejer1932b})}.
\end{equation}

These results have been extended to Jacobi-Gauss-Lobatto pointsystem (the roots of $(1-t^2)P_{n-2}^{(\alpha,\beta)}(t)=0$) and  Jacobi-Gauss-Radau pointsystem (the roots of $(1-t)P_{n-1}^{(\alpha,\beta)}(t)=0$ or $(1+t)P_{n-1}^{(\alpha,\beta)}(t)=0$) by V\'{e}rtesi \cite{Vertesi1979a,Vertesi1979b}: for all $k$ and $t\in [-1,1]$,
\begin{equation}
 v_k(t)\ge \min\{2-\alpha,2-\beta\}\mbox{\small\quad for $\{x_{k}\}\bigcup \{-1,1\}$ with  $1\le\alpha\le 2$ and $1\le\beta\le 2$,}
\end{equation}
\begin{equation}
 v_k(t)\ge \min\{2-\alpha,-\beta\}\mbox{\quad for $\{x_{k}\}\bigcup \{1\}$ with  $1\le\alpha\le 2$ and $-1<\beta\le 0$,}
\end{equation}
\begin{equation}
\quad\quad v_k(t)\ge \min\{-\alpha,2-\beta\}\mbox{\quad for $\{x_{k}\}\bigcup \{-1\}$ with  $-1<\alpha\le 0$ and $1\le\beta\le 2$}.
\end{equation}

\begin{proposition} (i) \cite{Fejer1932a,Szego} The Gauss-Jacobi pointsystem  is strongly normal if and only if  $\max\{\alpha,\beta\}<0$.

(ii) \cite{Vertesi1979a,Vertesi1979b} The Jacobi-Gauss-Lobatto pointsystem  is strongly  normal if and only if $1\le\alpha< 2$ and $1\le\beta< 2$.

(iii) \cite{Vertesi1979a,Vertesi1979b} The Jacobi-Gauss-Radau pointsystem including $x_1=1$  is strongly  normal if and only if  $1\le\alpha< 2$ and  $-1\le\beta< 0$, and the Jacobi-Gauss-Radau pointsystem including $x_n=-1$ is strongly  normal if and only if $-1<\alpha< 0$ and $1\le\beta< 2$.

\end{proposition}

It is worth noticing that if the pointsystem  is strongly normal, then it implies $v_i(t)\ge c>0$ for all $i=1,2,\ldots,n$ and $t\in [-1,1]$, and
\begin{equation}
1\equiv\sum_{i=1}^nh_i(t)=\sum_{i=1}^nv_i(t)\ell^2_i(t)\ge c\sum_{i=1}^n\ell^2_i(t)
\end{equation}
(see \cite{Fejer1932a}) and then
$$
\|\ell_i\|_{\infty}\le \frac{1}{\sqrt{c}},\quad i=1,2,\ldots,n.
$$

\begin{theorem}
Suppose $f(t)$ satisfies (2.9) and $\{x_j\}_{j=1}^n$ is a strongly normal pointsystem, then for $n\ge r+1$
\begin{equation}\label{error3-1}
\|E_n[f]\|_{\infty}\le \frac{\pi^r V_r}{\sqrt{c}(n-1)(n-2)\cdots(n-r)}.
\end{equation}
\end{theorem}

Following de la Vall\'{e}e Poussin \cite{Poussin1908}, the error bound indicates that $\|f-L_n[f]\|_{\infty}$ has the same asymptotic order as the estimate of $\|f-p_{n-1}^*\|_{\infty}$ for the interpolant at a strongly normal pointsystem for a functions of limited regularity with $V_r<\infty$ for some $r\ge 1$.

To check the error bounds in Theorem 4.2 numerically, we consider two limited regularity functions: $f(x)=|x|$ ($V_1<\infty$) and $f(x)=|x|^3$ ($V_3<\infty$).  All $(\alpha,\beta)$ are generated by ${\rm {\it rand}(1,2)}$ \footnote{${\rm {\it rand}(m,n)}$ returns an m-by-n matrix containing pseudo{\it rand}om values drawn from the standard uniform distribution on the open interval $(0,1)$.} except for $(\alpha,\beta)=(-0.5,-0.5)$, $(\alpha,\beta)=(0,0)$, $(\alpha,\beta)=(1,1)$ or $(\alpha,\beta)=(1.5,1.5)$. Particularly, we used ${\rm -{\it rand}(1,2)}$ in {\sc Figs}. 4.1-4.2 for strongly normal Gauss-Jacobi pointsystems, while ${\rm {\it rand}(1,2)}+1$ in {\sc Figs}. 4.3-4.4 for strongly normal Jacobi-Gauss-Lobatto pointsystems.
In {\sc Figs}. 4.5-4.6, we used $({\rm {\it rand}(1)+1}, -{\rm {\it rand}(1)})$ (1st row) and $({\rm -{\it rand}(1)}, {\rm {\it rand}(1)+1})$ (2nd row) for strongly normal Jacobi-Gauss-Radau pointsystems, respectively.

\begin{figure}[htbp]
\centerline{\includegraphics[height=6cm,width=14cm]{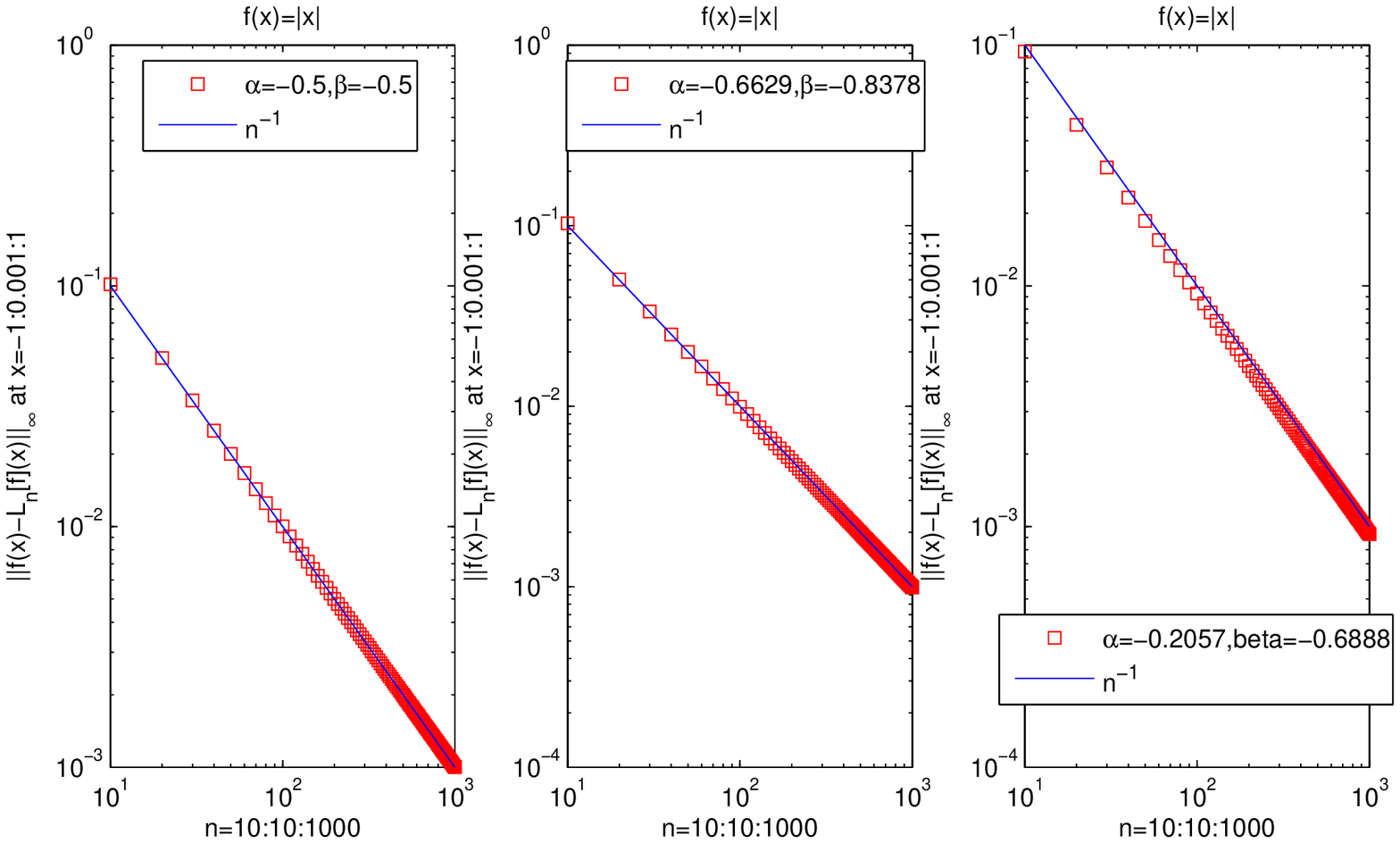}}
   \caption{$\max_{x=-1:0.001:1}|f(x)-L_n[f](x)|$ with $n=10:10:1000$ at the strongly normal Gauss-Jacobi pointsystems for $f(x)=|x|$, respectively.}
\end{figure}

\begin{figure}[htbp]
\centerline{\includegraphics[height=6cm,width=14cm]{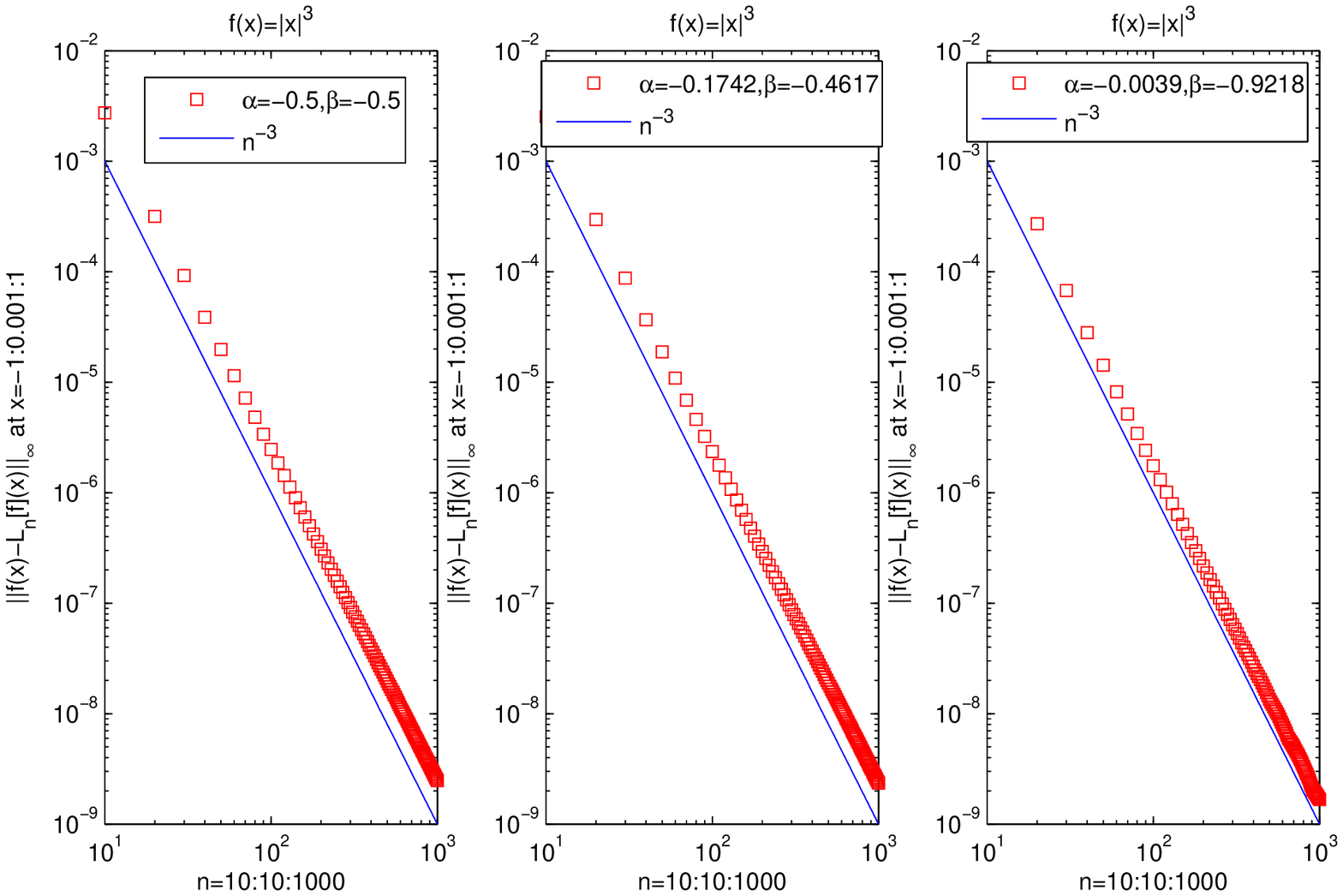}}
   \caption{$\max_{x=-1:0.001:1}|f(x)-L_n[f](x)|$ with $n=10:10:1000$ at the strongly normal Gauss-Jacobi pointsystems for $f(x)=|x|^3$, respectively.}
\end{figure}

\begin{figure}[htbp]
\centerline{\includegraphics[height=6cm,width=14cm]{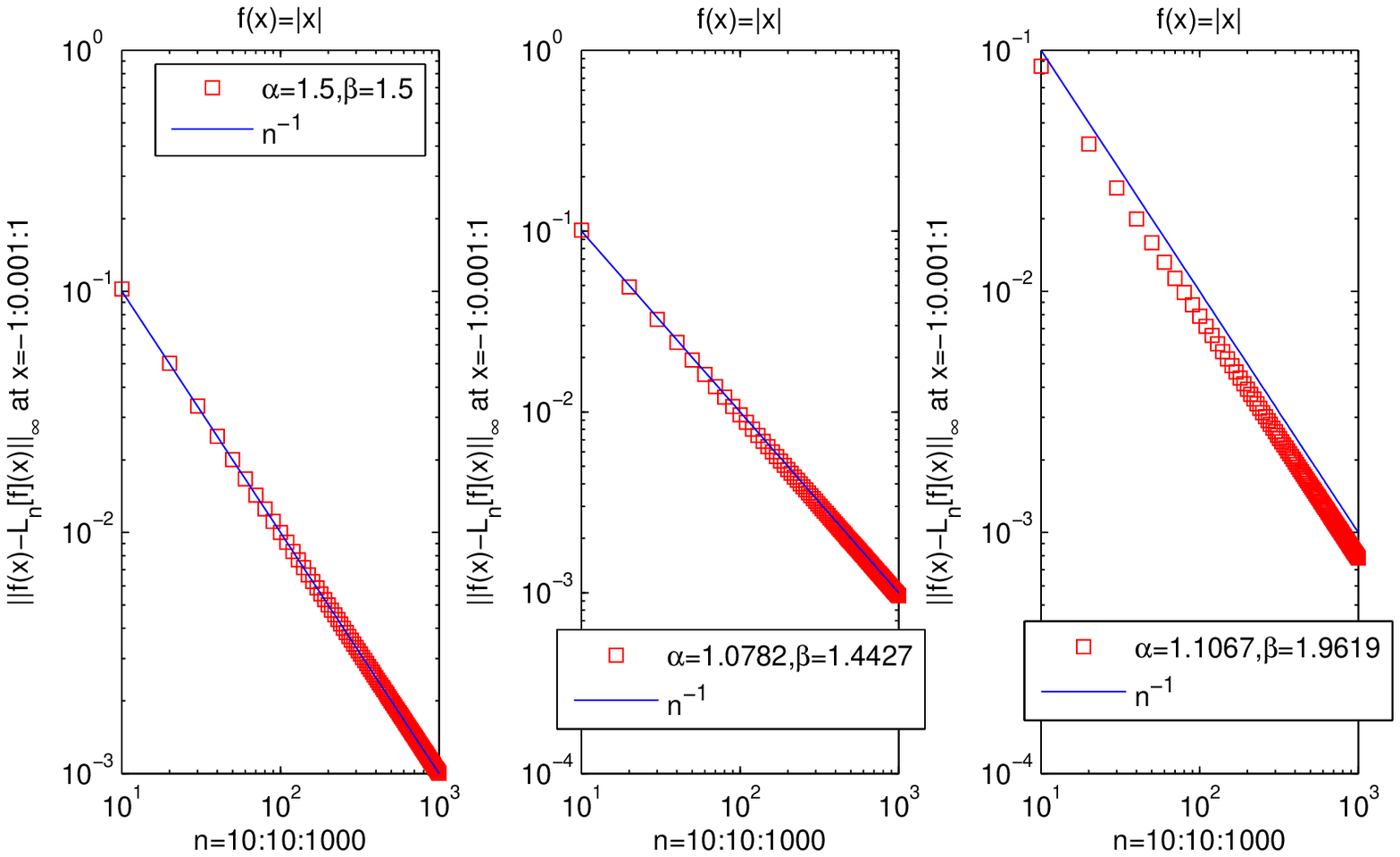}}
   \caption{$\max_{x=-1:0.001:1}|f(x)-L_n[f](x)|$ with $n=10:10:1000$ at the strongly normal Jacobi-Gauss-Lobatto pointsystems for $f(x)=|x|$, respectively.}
\end{figure}

\begin{figure}[htbp]
\centerline{\includegraphics[height=6cm,width=14cm]{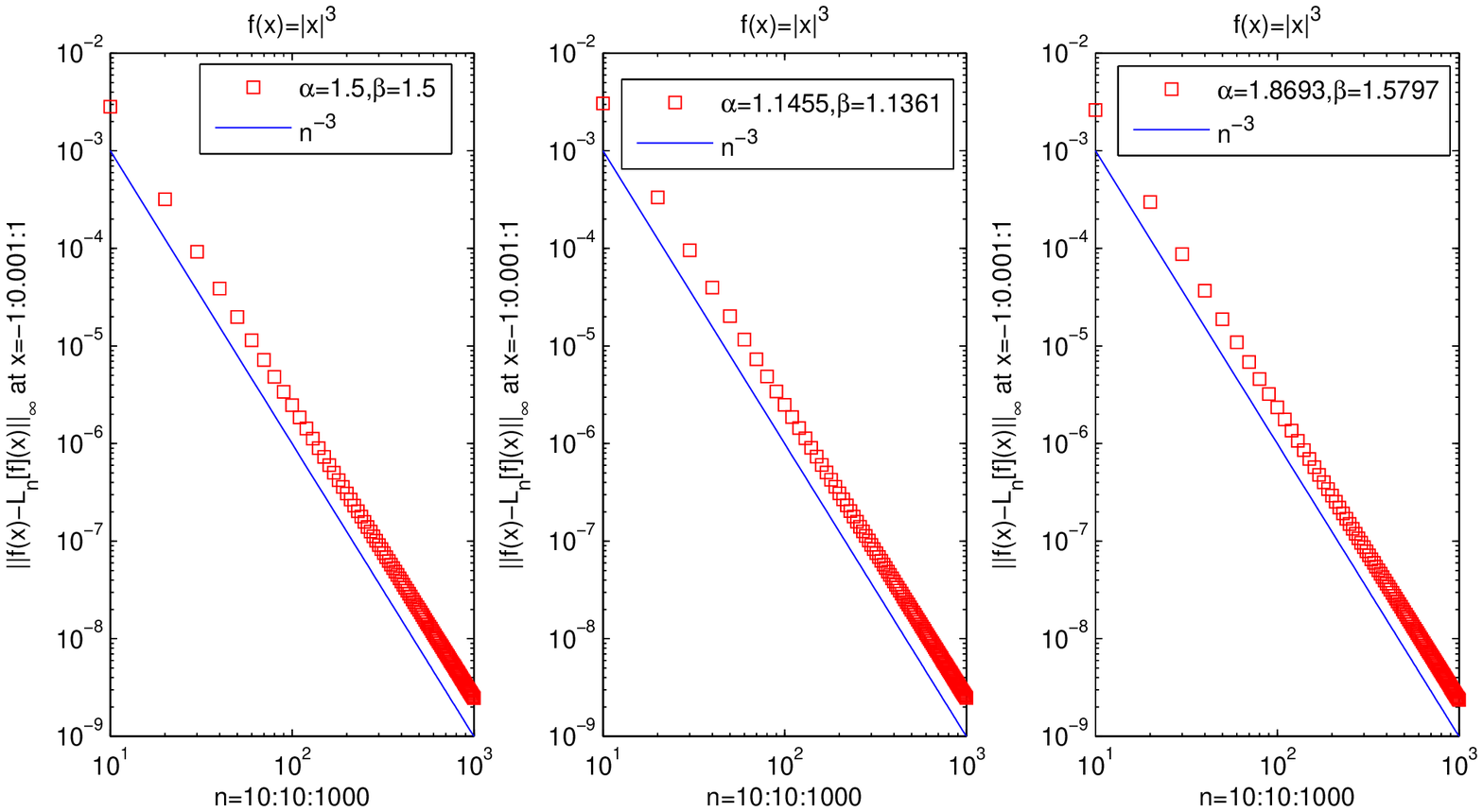}}
   \caption{$\max_{x=-1:0.001:1}|f(x)-L_n[f](x)|$ with $n=10:10:1000$ at the strongly normal Jacobi-Gauss-Lobatto pointsystems for $f(x)=|x|^3$, respectively.}
\end{figure}

\begin{figure}[htbp]
\centerline{\includegraphics[height=5cm,width=14cm]{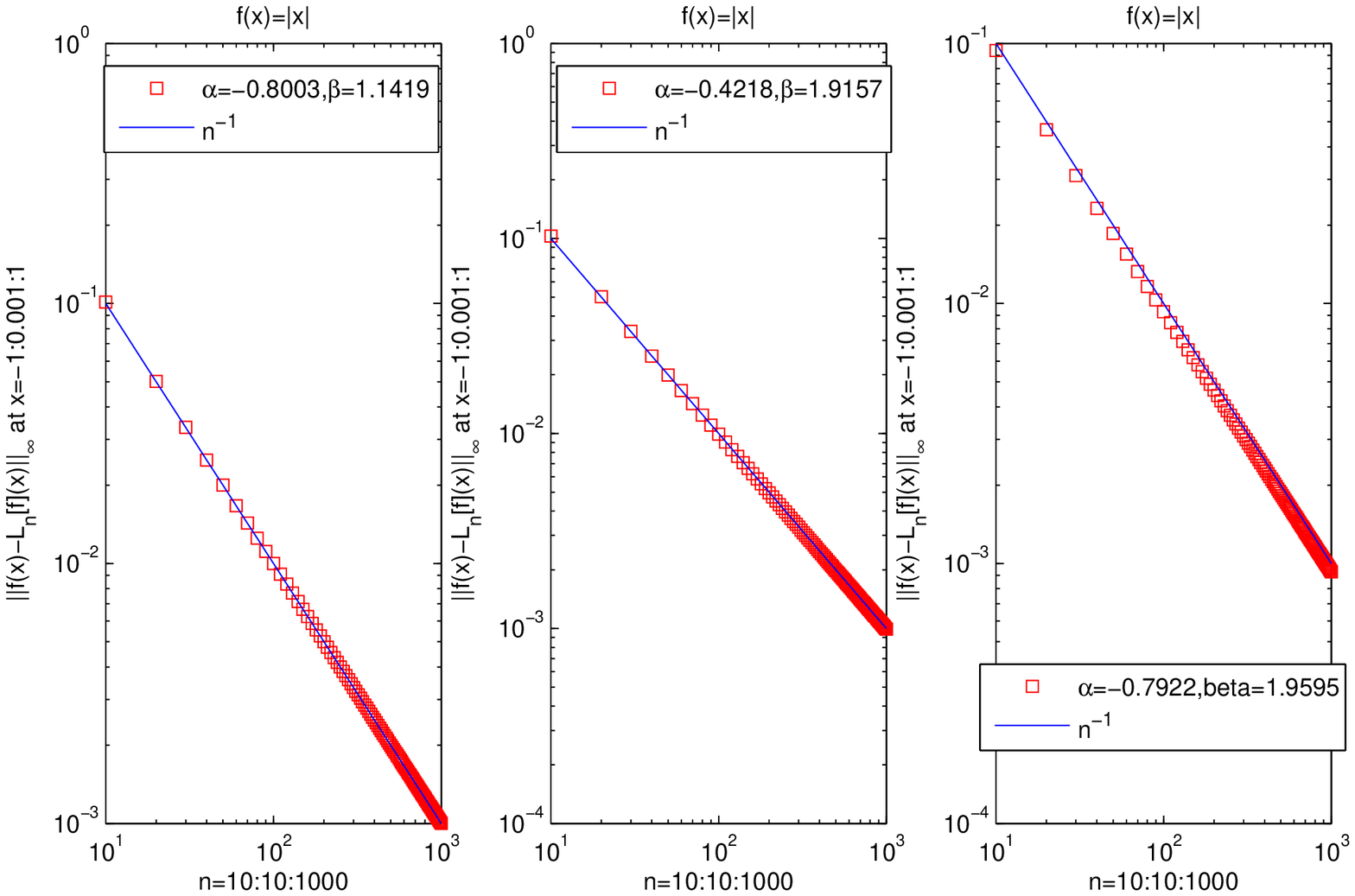}}
\centerline{\includegraphics[height=5cm,width=14cm]{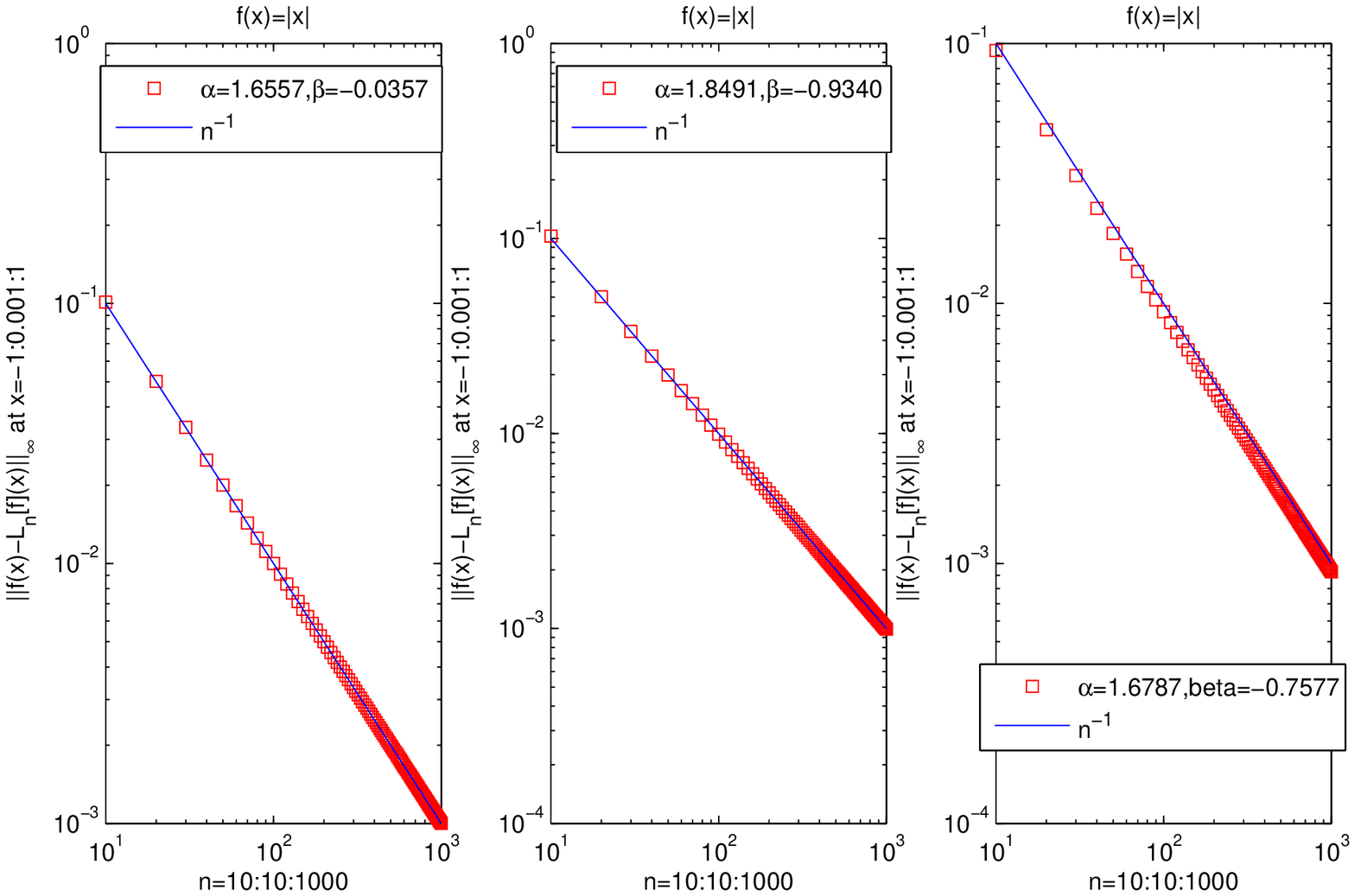}}
   \caption{$\max_{x=-1:0.001:1}|f(x)-L_n[f](x)|$ with $n=10:10:1000$ at the strongly normal Jacobi-Gauss-Radau pointsystems including $-1$ (1st row) and $1$ (2nd row) for $f(x)=|x|$, respectively.}
\end{figure}
\begin{figure}[htbp]
\centerline{\includegraphics[height=5cm,width=14cm]{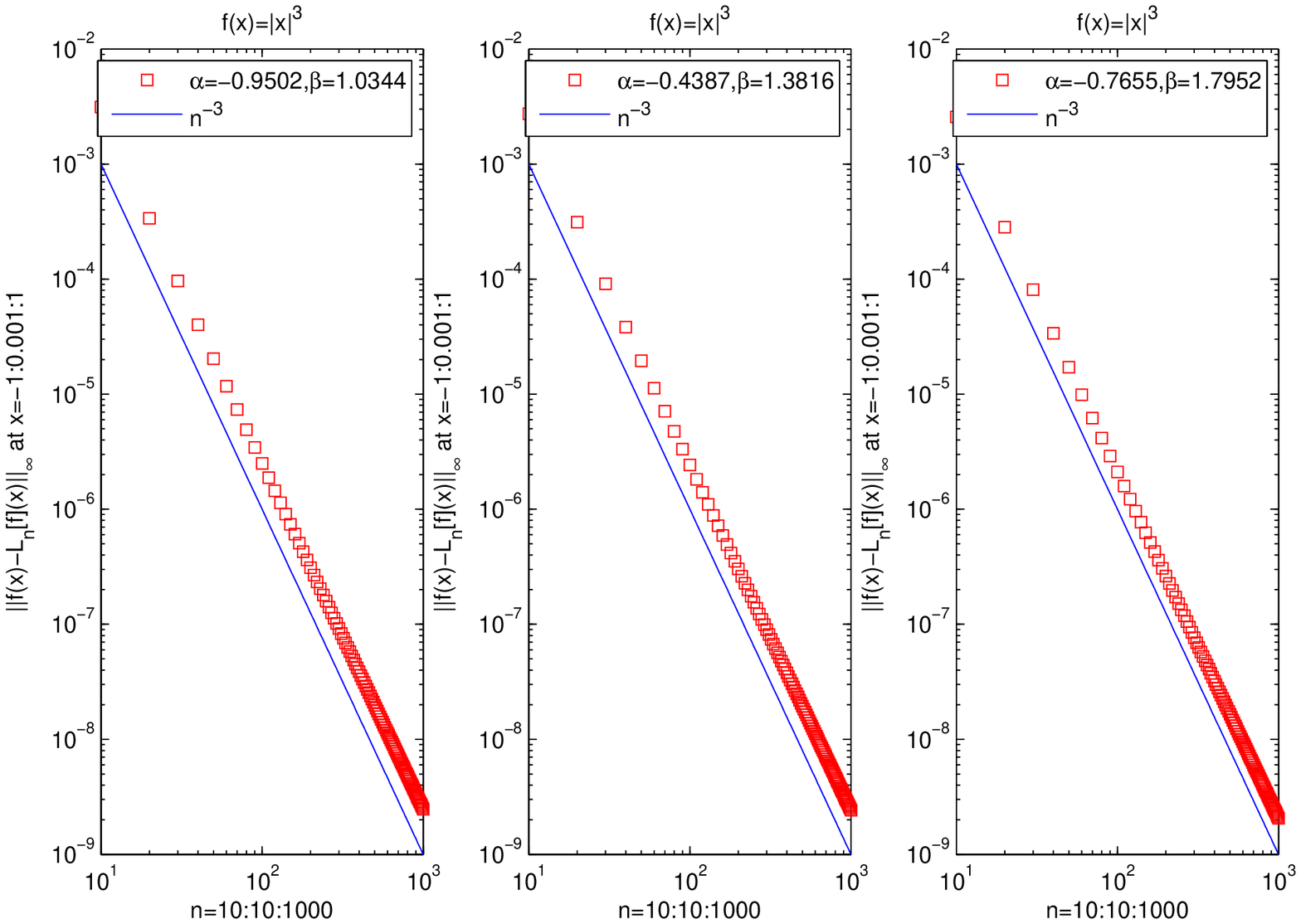}}\centerline{\includegraphics[height=5cm,width=14cm]{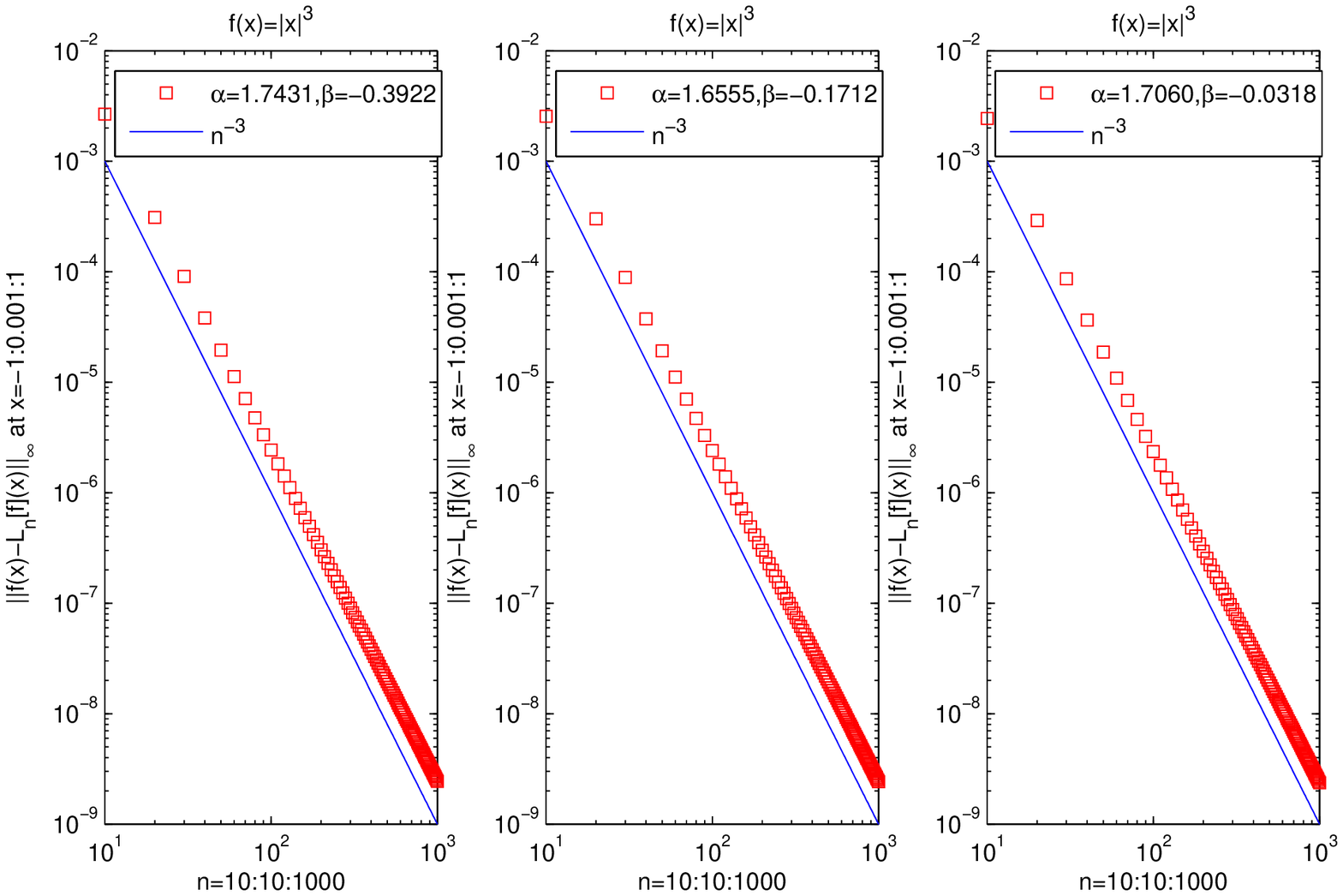}}
    \caption{$\max_{x=-1:0.001:1}|f(x)-L_n[f](x)|$ with $n=10:10:1000$ at the strongly normal Jacobi-Gauss-Radau pointsystems  including $-1$ (1st row) and $1$ (2nd row) for $f(x)=|x|^3$, respectively.}
\end{figure}

From {\sc Figs}. 4.1-4.6, we see that these convergence rates are in conformity to the estimates and attainable.


\subsection{General Gauss-Jacobi pointsystems}
In this subsection, we will consider convergence rates for general Gauss-Jacobi pointsystems, which includes the corresponding strongly normal pointsystems ($-1<\alpha,\beta<0$) as special cases.

Let $\{x_k\}_{k=1}^n$ be the roots of the Jacobi polynomial $P_n^{(\alpha,\beta)}(t)$ ($\alpha,\beta>-1$) and $x_k=\cos(\theta_k)$. Then from Szeg\"{o} \cite{Szego}, it follows
\begin{equation}
P_n^{(\alpha,\beta)}(t)=(-1)^nP_n^{(\beta,\alpha)}(-t)  \quad (\cite[(4.1.3)]{Szego})
\end{equation}
\begin{equation}
\max_{-1\le t\le 1}\big|P_n^{(\alpha,\beta)}(t)\big|=\left\{\begin{array}{ll}
\left(\begin{array}{c}n+q\\n\end{array}\right)\sim n^q,& q=\max\{\alpha,\beta\}\ge-\frac{1}{2}\\
\big|P_n^{(\alpha,\beta)}(t')\big|\sim n^{-\frac{1}{2}},& q=\max\{\alpha,\beta\}<-\frac{1}{2}
\end{array}\right.  \quad (\cite[(7.32.2)]{Szego})
\end{equation}
where $x'$ is one of the two maximum points, and for $t=\cos(\theta)$ and  any fixed  constant $c$ with $0<c<1$,
\begin{equation}
P_{n}^{(\alpha,\beta)}(\cos(\theta))=\left\{\begin{array}{ll}
O\left(n^{\alpha}\right),& 0\le \theta\le cn^{-1}\\
\theta^{-\alpha-\frac{1}{2}}O\left(n^{-\frac{1}{2}}\right),& cn^{-1}\le \theta\le\frac{\pi}{2}\end{array}\right.\quad (\mbox{\cite[Theorem 7.32.2]{Szego}}),
\end{equation}

\begin{equation}
\theta_k=n^{-1}\left[k\pi+O(1)\right] \quad (\cite[(8.9.1)]{Szego}),
\end{equation}
\begin{equation}
{\big|P_n^{(\alpha,\beta)}}'(\cos(\theta_k))\big|\sim k^{-\alpha-\frac{3}{2}}n^{\alpha+2}, \quad 0<\theta_k\le\frac{\pi}{2} \quad (\cite[(8.9.2)]{Szego}).
\end{equation}
Moreover, expression (4.17) can be extended to
\begin{equation}\quad
{\big|P_n^{(\alpha,\beta)}}'(\cos(\theta_k))\big|\sim k^{-\alpha-\frac{3}{2}}n^{\alpha+2}, \, 0<\theta_k\le c_1\pi
\end{equation}
for any fixed $c_1$ with $0<c_1<1$ (\cite[(4.6)]{Vertesi1983}).

Based on these identities, the estimates on  ${\displaystyle \ell_k(t)=\frac{P_n^{(\alpha,\beta)}(t)}{{P_n^{(\alpha,\beta)}}'(x_k)(t-x_k)}}$ have been extensively studied in Kelzon \cite{Kelzon,Kelzon2}, V\'{e}rtesi \cite{Vertesi1980,Vertesi1983b}, Sun \cite{Sun}, Prestin \cite{Prestin}, Kvernadze \cite{Kvernadze}, Vecchia et al. \cite{Vecchia}, etc.

\begin{lemma} \cite{Sun} (also see \cite{Kvernadze})
For $t\in [-1,1]$, let $x_m$ be the root of the Jacobi polynomial $P_n^{(\alpha,\beta)}$ which is closest to $t$. Then we have
\begin{equation}\label{Sun}
  \ell_k(t)=\left\{\begin{array}{ll}O\left(|k-m|^{-1}+|k-m|^{\gamma-\frac{1}{2}}\right),&k\not=m\\
  O(1)&k=m\end{array}\right.,\quad  \gamma=\max\{\alpha,\beta\},
\end{equation}
for $k=1,2,\ldots,n$.
\end{lemma}
\begin{proof}
In \cite{Sun}, the proof of Lemma 4.3 is given only for $0\le \theta_k\le \frac{\pi}{2}$ or $k=m$. That proof can be readily extended to $0\le \theta_k\le \frac{2\pi}{3}$ due to (4.18). We complement the proof for $\frac{2\pi}{3}<\theta_k< \pi$ and $k\not=m$ next.

From (4.13) and (4.18), we see that
\begin{equation}
{\big|P_n^{(\alpha,\beta)}}'(\cos(\theta_k))\big|\sim (n-k+1)^{-\beta-\frac{3}{2}}n^{\beta+2}, \quad \frac{2\pi}{3}<\theta_k< \pi \quad (\cite[(9)]{Prestin}).
\end{equation}
Then for $0\le t=\cos (\theta)\le 1$ with $0\le \theta\le cn^{-1}$ and $\frac{2\pi}{3}<\theta_k< \pi$, it follows by (4.15) and (4.20) that
$$
\ell_k(t)=O\left(\frac{n^{\alpha}}{(n-k+1)^{-\beta-\frac{3}{2}}n^{\beta+2}}\right)=O\left(\frac{(n-k+1)^{\beta+\frac{3}{2}}}{n^{\beta+2-\alpha}}\right)
=O\left(n^{\alpha-\frac{1}{2}}\right).
$$
While for $cn^{-1}\le \theta\le \frac{\pi}{2}$ and $\frac{2\pi}{3}<\theta_k< \pi$, it follows by (4.15)-(4.18) and (4.20) that
$$
\ell_k(t)=O\left(\frac{(m\pi/n)^{-\alpha-\frac{1}{2}}n^{-\frac{1}{2}}}{(n-k+1)^{-\beta-\frac{3}{2}}n^{\beta+2}}\right)
=O\left(\frac{1}{m^{\frac{1}{2}+\alpha}n^{\frac{1}{2}-\alpha}}\right)=\left\{\begin{array}{ll}
O\left(n^{\alpha-\frac{1}{2}}\right),&\alpha>-\frac{1}{2}\\
O\left(n^{-1}\right),&-1<\alpha\le -\frac{1}{2}.\end{array}\right.
$$
Thus for $0\le t\le 1$, we have $\ell_k(t)=O\left(n^{-1}+n^{\alpha-\frac{1}{2}}\right)$ for $k\not=m$, which leads to the desired result due to that $k-m\sim n$ in the case $\frac{2\pi}{3}<\theta_k< \pi$.

Similarly, by (4.13) together with the above analysis, we get for $-1\le t\le 0$ that $$\ell_k(t)=O\left(|k-m|^{-1}+|k-m|^{\beta-\frac{1}{2}}\right),\quad k\not=m.$$
These together lead to the desired result (4.19) for $k\not=m$.
\end{proof}

\begin{theorem}
Suppose $f(t)$ satisfies (2.9) and $\{x_j\}_{j=1}^n$  are the roots of the Jacobi polynomial $P_n^{(\alpha,\beta)}(t)$, then for $n\ge r+1$
\begin{equation}\label{error4-1}
\|E_n[f]\|_{\infty}=O\left(n^{-r+\max\left\{0,\gamma-\frac{1}{2}\right\}}\right),\quad \gamma=\max\{\alpha,\beta\}
\end{equation}
\end{theorem}
\begin{proof}
From Lemma 4.3, we see that $\max_{1\le j\le n}\|\ell_j\|_{\infty}=O\left(n^{\max\left\{0,\gamma-\frac{1}{2}\right\}}\right)$, which together with Theorem 3.2 yields the desired result.
\end{proof}

{\sc Remark 2}.  Theorem 4.4 implies that $\|f-L_n[f]\|_{\infty}$ has the same asymptotic order as $\|f-p_{n-1}^*\|_{\infty}$ \cite{Poussin1908}
at the roots of the Jacobi polynomial $P_n^{(\alpha,\beta)}(t)$ for $-1<\alpha,\beta\le \frac{1}{2}$. Then the interpolations at the $n$-point Gauss-Legendre points and at the $n$-point Chebyshev points of first kind or second kind have essentially the same accuracy. All of them can achieve the optimal convergence rate  $O(\|f-p_{n-1}^*\|_{\infty})$. Consequently, the corresponding quadrature Gauss, Clenshaw-Curtis and Fej\'{e}r first rule have essentially the same accuracy \cite{Xiang0}.

Here, we used {\sc Figs}. 4.7-4.8 to illustrate the convergence rates for general Gauss-Jacobi pointsystems, where $(\alpha,\beta)$ are obtained by ${\rm {\it {\it rand}}(1,2)}$ (1st row) and $m{\rm {\it {\it rand}}(1,2)}$ with $m\|{\rm {\it {\it rand}}}(1,2)\|_{\infty}$ $>m-1$ for $m=2,3,4$ (2nd row), respectively. From these figures, we see that the convergence rates are attainable too, which are in accordance with the estimates. Then the convergence rates at the Gauss-Jacobi pointsystems are optimal.

\begin{figure}[htbp]
\centerline{\includegraphics[height=8cm,width=14cm]{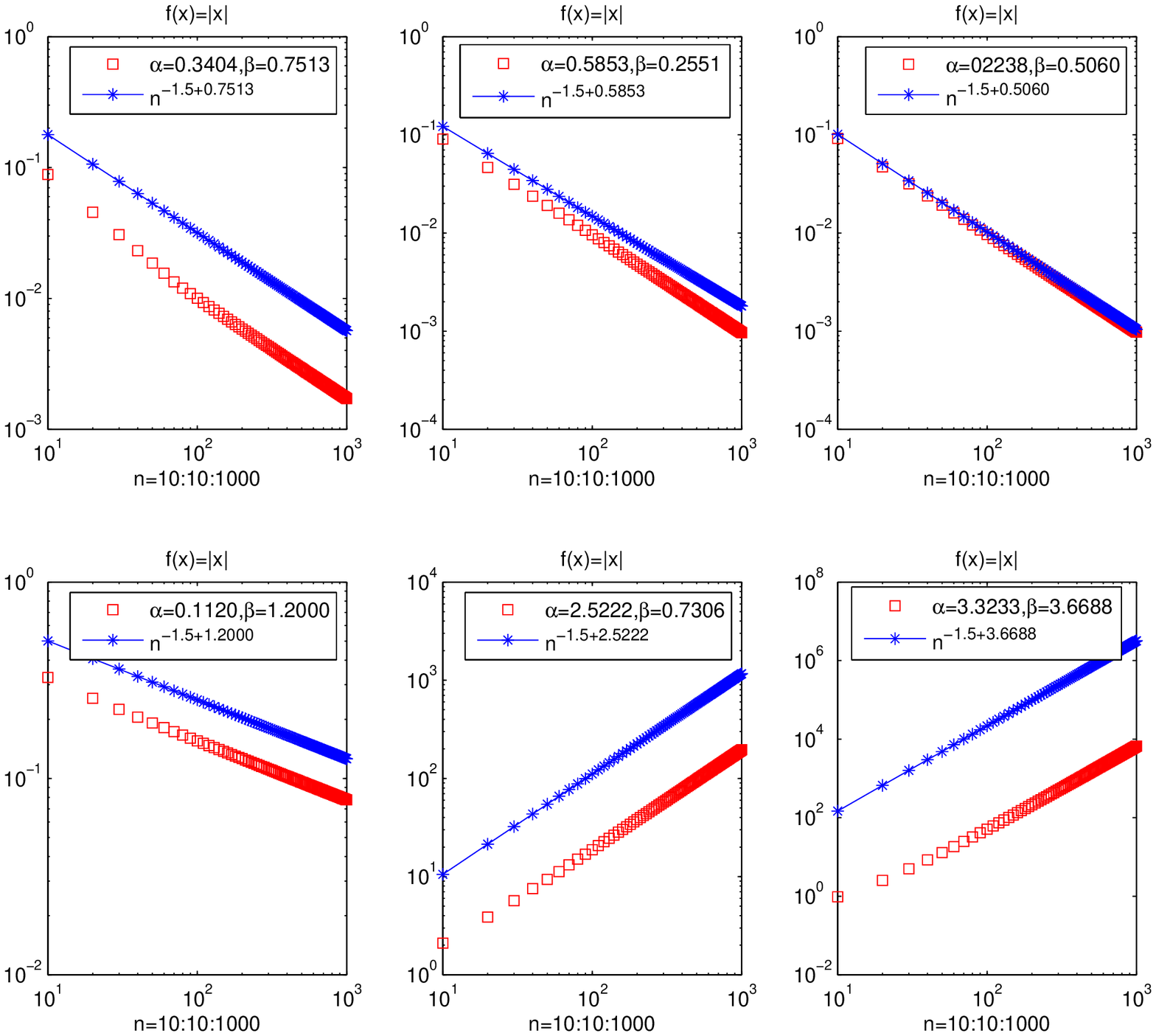}}
\vspace{-.6cm}
   \caption{$\max_{x=-1:0.001:1}|f(x)-L_n[f](x)|$ with $n=10:10:1000$ at the Gauss-Jacobi pointsystems for $f(x)=|x|$, respectively.}
\end{figure}

\begin{figure}[htbp]
\centerline{\includegraphics[height=8cm,width=14cm]{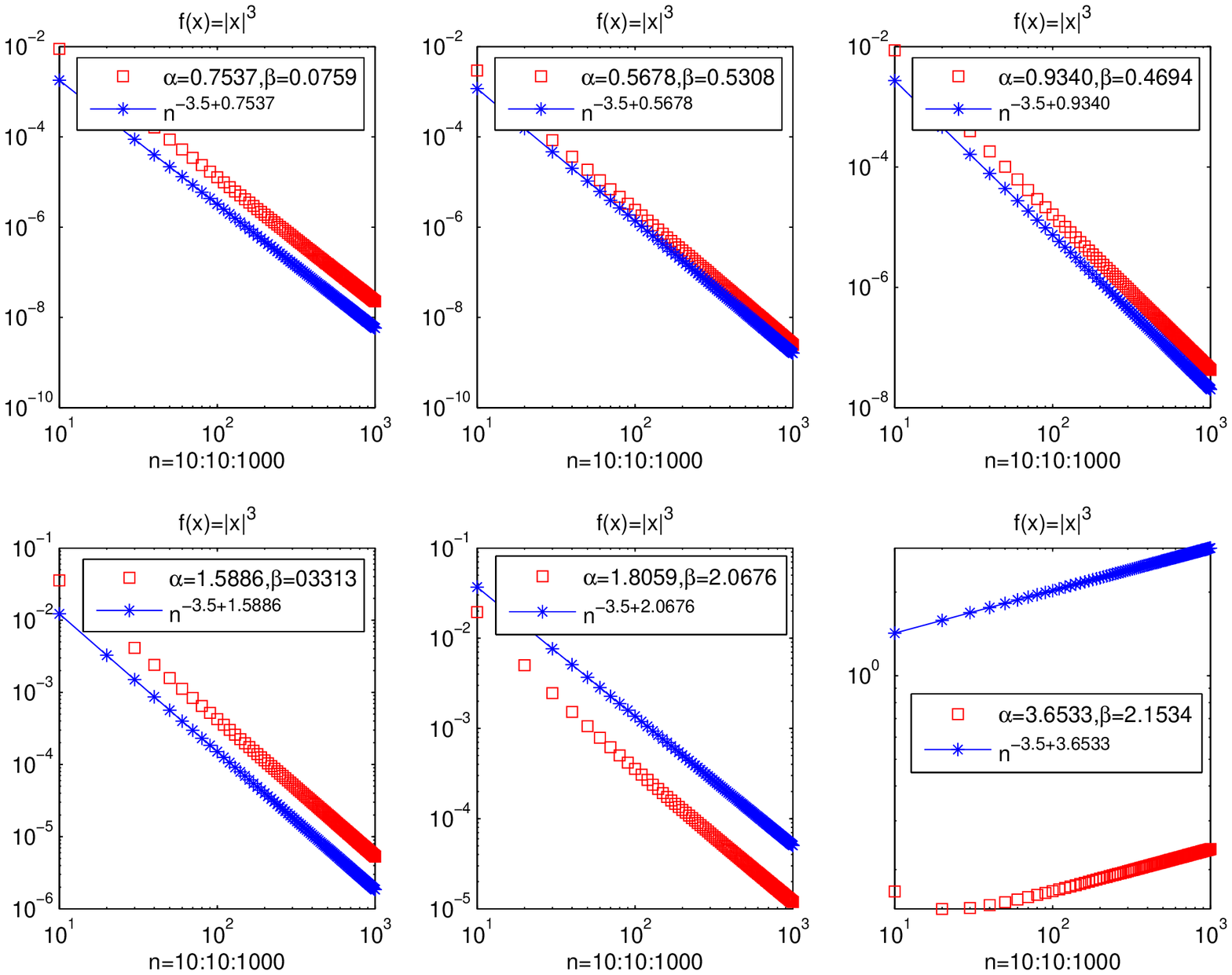}}
   \caption{$\max_{x=-1:0.001:1}|f(x)-L_n[f](x)|$ with $n=10:10:1000$ at the Gauss-Jacobi pointsystems for $f(x)=|x|^3$, respectively.}
\end{figure}

{\sc Remark 3}. It is of particular relevance from {\sc Figs}. 4.7-4.8  in the cases that  the polynomial interpolations are divergent  if $r-\max\left\{0,\gamma-\frac{1}{2}\right\}\le 0$, the divergence rate is also controlled by the order $O\left(n^{-r+\max\left\{0,\gamma-\frac{1}{2}\right\}}\right)$.


\subsection{General Jacobi-Gauss-Lobatto pointsystems}
Let
\begin{equation}\label{LobattoGaussJacobi}
-1=x_{n+1}<x_n<x_{n-1}<\cdots<x_2<x_1<x_0=1
\end{equation}
be the roots of $(1-t^2)P_{n}^{(\alpha,\beta)}(t)=0$ ($\alpha,\beta>-1$), $x_k=\cos(\theta_k)$ and
$$
\omega(t)=(t-x_0)(t-x_1)\cdots(t-x_n)(t-x_{n+1}), \quad \ell_k(t)=\frac{\omega(t)}{(t-x_k)\omega'(x_k)}.$$
Then
\begin{equation}
\ell_0(t)=\frac{1+t}{2}\cdot\frac{P_{n}^{(\alpha,\beta)}(t)}{P_{n}^{(\alpha,\beta)}(1)},\quad \ell_{n+1}(t)=\frac{1-t}{2}\cdot\frac{P_{n}^{(\alpha,\beta)}(t)}{P_{n}^{(\alpha,\beta)}(-1)},
\end{equation}
and
\begin{equation}
{\displaystyle\ell_k(t)=\frac{(1-t^2)P_{n}^{(\alpha,\beta)}(t)}{(t-x_k)(1-x_k^2){P_{n}^{(\alpha,\beta)}}'(x_k)}},\quad k=1,2,\ldots,n.
\end{equation}
In the next, we shall concentrate on estimates of $\ell_k(t)$ for $k=0,1,2,\ldots,n+1$.

\begin{itemize}
\item On the estimate of $\ell_0(t)$: (i) In the case $0\le t\le 1$, setting $t=\cos\theta$ for $0\le \theta\le \frac{\pi}{2}$, and using $$P_{n}^{(\alpha,\beta)}(1)=\left(\begin{array}{c}n+\alpha\\n\end{array}\right)\sim n^{\alpha}\,\,\, (\cite[(4.1.1),(7.32.2)]{Szego}),$$we find that  from (4.15) and (4.23) for $0\le \theta\le \frac{\pi}{2}$,
$$
\ell_0(t)=\left\{\begin{array}{ll}
O(1),& 0\le \theta\le cn^{-1}\\
O\left(\theta^{-\alpha-\frac{1}{2}}n^{-\frac{1}{2}}n^{-\alpha}\right)=O\left((n\theta)^{-\alpha-\frac{1}{2}}\right)=
O\left(n^{-\min\{0,\alpha+\frac{1}{2}\}}\right),& cn^{-1}\le \theta\le\frac{\pi}{2}\end{array}.\right.
$$
(ii) In the case $-1\le t\le 0$, letting $t=-\cos(\theta)$ for $0\le \theta\le \frac{\pi}{2}$ and applying $P_{n}^{(\alpha,\beta)}(-\cos(\theta))=(-1)^nP_{n}^{(\beta,\alpha)}(\cos(\theta))$ and
$1-\cos(\theta)=2\sin^2\left(\frac{\theta}{2}\right)$ and $\frac{2}{\pi}(\theta)\le\sin(\theta)\le \theta$, together with (4.15) and (4.23), we have
$$\begin{array}{lll}
\ell_0(t)&=&O\left(\theta^2P_{n}^{(\beta,\alpha)}(\cos(\theta))n^{-\alpha}\right)\\
&=&
\left\{\begin{array}{ll}
O\left(\frac{1}{n^{2+\alpha-\beta}}\right),& 0\le \theta\le cn^{-1}\\
O\left(\theta^{-\beta+\frac{3}{2}}n^{-\frac{1}{2}}n^{-\alpha}\right)=O\left(n^{-\min\{2+\alpha-\beta,\alpha+\frac{1}{2}\}}\right),& cn^{-1}\le \theta\le\frac{\pi}{2}\end{array}.\right.\end{array}
$$
These together yield
\begin{equation}
\|\ell_0\|_{\infty}=O\left(\frac{1}{n^{\min\left\{0,2+\alpha-\beta,\alpha+\frac{1}{2}\right\}}}\right).
\end{equation}

\item Similarly, we have
\begin{equation}
\|\ell_{n+1}\|_{\infty}=O\left(\frac{1}{n^{\min\left\{0,2+\beta-\alpha,\beta+\frac{1}{2}\right\}}}\right).
\end{equation}

\item  For $k=1,2,\ldots,n$, let $x_m$ be the nearest to $t\in [0,1]$ and $t=\cos(\theta)$.  From (4.24), we have for $k\not=m$ that
\begin{equation}\quad
\ell_k(t)=\frac{\sin^2\theta P_{n}^{(\alpha,\beta)}(\cos\theta)}{(\cos\theta-\cos\theta_k)\sin^2\theta_k{P_{n}^{(\alpha,\beta)}}'(\cos\theta_k)}=\frac{-\sin^2\theta P_{n}^{(\alpha,\beta)}(\cos\theta)}{2\sin\left(\frac{\theta-\theta_k}{2}\right)\sin\left(\frac{\theta+\theta_k}{2}\right)\sin^2\theta_k
{P_{n}^{(\alpha,\beta)}}'(\cos\theta_k)}.
\end{equation}
\textbf{In the case $0\le \theta\le cn^{-1}$ and $0\le \theta_k\le \frac{2\pi}{3}$}: From (4.15)-(4.18), it follows
\begin{equation}\quad\quad
\ell_k(\cos\theta)=O\left(\frac{n^{-2}n^{\alpha}}{|k-m||k+m|n^{-2}k^2n^{-2}k^{-\alpha-\frac{3}{2}}n^{\alpha+2}}\right)
=O\left(\frac{k^{\alpha-\frac{1}{2}}}{|k-m||k+m|}\right).
\end{equation}
Define
$$
h_1(u)=\frac{u^{\alpha-\frac{1}{2}}}{u^2-m^2} \mbox{\quad for $m+1\le u\le n$};\quad h_2(u)=-\frac{u^{\alpha-\frac{1}{2}}}{u^2-m^2} \mbox{\quad for $1\le u\le m-1$}.
$$
Then by an elementary proof and noting that $m\le c_1n$ for $0<c_1<1$, we get
$$
\max_{m+1\le u\le n} h_1(u)=\left\{\begin{array}{ll}
h_1(m+1)=O\left(m^{\alpha-\frac{3}{2}}\right),&-1<\alpha\le \frac{5}{2}\\
\max\left\{h_1(m+1),h_1(n)\right\}=O\left(\max\left\{m^{\alpha-\frac{3}{2}},n^{\alpha-\frac{5}{2}}\right\}\right),&\alpha>\frac{5}{2}\end{array}\right.
$$
and
$$
\max_{1\le u\le m-1} h_2(u)=\left\{\begin{array}{ll}
\max\left\{h_2(1),h_2(m-1)\right\}=O\left(\max\left\{m^{-2}, m^{\alpha-\frac{5}{2}}\right\}\right),&-1<\alpha\le \frac{1}{2}\\
h_2(m-1)=O\left(m^{\alpha-\frac{3}{2}}\right),&\alpha> \frac{1}{2}\end{array}\right.,
$$
which, together with $m\sim 1$ under the assumption, establishes that
\begin{equation}
\ell_k(\cos\theta)=\left\{\begin{array}{ll}
O\left(1\right),&-1<\alpha\le \frac{5}{2}\\
O\left(n^{\alpha-\frac{5}{2}}\right),&\alpha> \frac{5}{2}\end{array}.\right.
\end{equation}
\textbf{In the case $0\le \theta\le cn^{-1}$ and $\frac{2\pi}{3}< \theta_k<\pi$}: Similarly, from (4.15) and (4.20) we have
\begin{equation}\begin{array}{lll}
\ell_k(\cos\theta)&=&O\left(\frac{n^{-2}n^{\alpha}}{(n-k+1)^2n^{-2}(n-k+1)^{-\beta-\frac{3}{2}}n^{\beta+2}}\right)\\
&=&O\left(\frac{(n-k+1)^{\beta-\frac{1}{2}}}{n^{2+\beta-\alpha}}\right)\\
&=&O\left(n^{-\min\{2+\beta-\alpha,\frac{5}{2}-\alpha\}}\right).\end{array}
\end{equation}

\textbf{In the case $cn^{-1}\le \theta\le \frac{\pi}{2}$  and $0\le \theta_k\le \frac{2\pi}{3}$}: By (4.27), together with (4.15)-(4.18) , we obtain
\begin{equation}
\ell_k(\cos\theta)=O\left(\frac{m^{\frac{3}{2}-\alpha}k^{\alpha-\frac{1}{2}}}{|k-m||k+m|}\right)
=m^{\frac{3}{2}-\alpha}O\left(\frac{k^{\alpha-\frac{1}{2}}}{|k-m||k+m|}\right)
\end{equation}
which  establishes that by applying the estimates to $h_1(u)$ and $h_2(u)$
\begin{equation}
\ell_k(\cos\theta)=\left\{\begin{array}{ll}
O\left(n^{-\alpha-\frac{1}{2}}\right),&-1<\alpha< -\frac{1}{2}\\
O\left(1\right),&-\frac{1}{2}\le\alpha\le \frac{5}{2}\\
O\left(n^{\alpha-\frac{5}{2}}\right),&\alpha> \frac{5}{2}\end{array}.\right.
\end{equation}
\textbf{In the case $cn^{-1}\le \theta\le \frac{\pi}{2}$  and $\frac{2\pi}{3}< \theta_k<\pi $}: From (4.15)-(4.18), (4.20) and (4.27), we find that
\begin{equation}\quad
\ell_k(\cos\theta)=m^{\frac{3}{2}-\alpha}O\left(\frac{(n-k+1)^{\beta-\frac{1}{2}}}{n^{2+\beta-\alpha}}\right)
=\left\{\begin{array}{ll}
O\left(n^{-\min\{0,\beta+\frac{1}{2}\}}\right),&-1<\alpha\le \frac{3}{2}\\
O\left(n^{-\min\{2+\beta-\alpha,\frac{5}{2}-\alpha\}}\right),&\alpha>\frac{3}{2}\end{array}.\right.
\end{equation}
Thus for $t\in [0,1]$, we get
\begin{equation}
\|\ell_k\|_{\infty}=\left\{\begin{array}{ll}
O\left(n^{-\min\{0,\beta+\frac{1}{2}, \alpha+\frac{1}{2}\}}\right),&-1<\alpha\le \frac{3}{2}\\
O\left(n^{-\min\{0,2+\beta-\alpha,\frac{5}{2}-\alpha\}}\right),&\alpha>\frac{3}{2}\end{array}.\right.
\end{equation}

For $t\in [-1,0]$, by $P_{n}^{(\beta,\alpha)}(-t)=(-1)^nP_{n}^{(\alpha,\beta)}(t)$, setting $t=-\cos\theta$ and $y_k=-x_{n-k+1}=\cos{\overline{\theta}_k}$ for $k=1,2,\ldots,n$, we see that $y_k$ are the roots of $P_{n}^{(\alpha,\beta)}(-t)=(-1)^nP_{n}^{(\beta,\alpha)}(t)$, then (4.24) can be represented as for $ k=1,2,\ldots,n$
$$
{\displaystyle\ell_{n-k+1}(t)=\frac{(1-\cos^2\theta)(-1)^nP_{n}^{(\beta,\alpha)}(\cos\theta)}{-(\cos\theta-y_{k})(1-y_{k}^2){P_{n}^{(\alpha,\beta)}}'(-y_{k})}
=\frac{\sin^2\theta P_{n}^{(\beta,\alpha)}(\cos\theta)}{(\cos\theta-y_{k})(1-y_{k}^2){P_{n}^{(\beta,\alpha)}}'(y_{k})}}
$$
followed
$$
\left[P_{n}^{(\alpha,\beta)}(t)\right]'=\frac{1}{2}(n+\alpha+\beta+1)P_{n-1}^{(\alpha+1,\beta+1)}(t)\,\,\cite[(4.21.7)]{Szego}.
$$
Similarly, we get
that
\begin{equation}
\|\ell_k\|_{\infty}=\left\{\begin{array}{ll}
O\left(n^{-\min\{0,\alpha+\frac{1}{2}, \beta+\frac{1}{2}\}}\right),&-1<\beta\le \frac{3}{2}\\
O\left(n^{-\min\{0,2+\alpha-\beta,\frac{5}{2}-\beta\}}\right),&\beta>\frac{3}{2}\end{array},\right.
\end{equation}
which together with (4.34) leads to that for $t\in [-1,1]$
\begin{equation}\quad\quad
\|\ell_k\|_{\infty}=\left\{\begin{array}{ll}
O\left(n^{-\min\{0,\alpha+\frac{1}{2},\beta+\frac{1}{2}\}}\right),&-1<\alpha,\beta\le \frac{3}{2}\\
O\left(n^{-\min\{0,\alpha+\frac{1}{2},2+\alpha-\beta,\frac{5}{2}-\beta\}}\right),&-1<\alpha\le \frac{3}{2},\beta> \frac{3}{2}\\
O\left(n^{-\min\{0,\beta+\frac{1}{2},2+\beta-\alpha,\frac{5}{2}-\alpha\}}\right),&\alpha>\frac{3}{2}, -1<\beta\le\frac{3}{2}\\
O\left(n^{-\min\{0,2+\alpha-\beta,2+\beta-\alpha,\frac{5}{2}-\alpha,\frac{5}{2}-\beta\}}\right),&\alpha,\beta>\frac{3}{2}\end{array}.\right.
\end{equation}
\end{itemize}

\begin{theorem}
Suppose $f(t)$ satisfies (2.9) and $\{x_j\}_{j=0}^{n+1}$  are the roots of $(1-t^2)P_n^{(\alpha,\beta)}(t)$ , then for $n\ge r+1$
\begin{equation}
E_n[f]=n^{-r}\cdot\left\{\begin{array}{ll}
O\left(n^{-\min\{0,\alpha+\frac{1}{2},\beta+\frac{1}{2}\}}\right),&-1<\alpha,\beta\le \frac{3}{2}\\
O\left(n^{-\min\{0,\alpha+\frac{1}{2},2+\alpha-\beta,\frac{5}{2}-\beta\}}\right),&-1<\alpha\le \frac{3}{2},\beta> \frac{3}{2}\\
O\left(n^{-\min\{0,\beta+\frac{1}{2},2+\beta-\alpha,\frac{5}{2}-\alpha\}}\right),&\alpha>\frac{3}{2}, -1<\beta\le\frac{3}{2}\\
O\left(n^{-\min\{0,2+\alpha-\beta,2+\beta-\alpha,\frac{5}{2}-\alpha,\frac{5}{2}-\beta\}}\right),&\alpha,\beta>\frac{3}{2}\end{array}.\right.\end{equation}
Particularly, we have for $(\alpha,\beta)\in S$
\begin{equation}
\|E_n[f]\|_{\infty}=O\left(\frac{1}{n^{r}}\right),
\end{equation}
where $S:= \left[-\frac{1}{2},\frac{5}{2}\right]\times \left[-\frac{1}{2},\frac{5}{2}\right]-\{(\alpha,\beta): -\frac{1}{2}\le \alpha\le \frac{1}{2},\,  2+\alpha<\beta\le \frac{5}{2}\}\cup \{(\alpha,\beta): \frac{3}{2}\le \alpha\le \frac{5}{2},\, -\frac{1}{2}\le\beta<2-\alpha\}$.

\end{theorem}

{\sc Remark 3}.  Theorem 4.5  implies $\|f-L_n[f]\|_{\infty}$ has the same asymptotic order as $\|f-p_{n-1}^*\|_{\infty}$ \cite{Poussin1908} at the roots of the Jacobi polynomial $(1-t^2)P_n^{(\alpha,\beta)}(t)$ for $(\alpha,\beta)\in S$, which includes the corresponding strongly normal pointsystems as special cases.

{\sc Figs}. 4.9-4.10 show the convergence rates for $f(x)=|x|$ or $f(x)=|x|^3$ at the Jacobi-Gauss-Lobatto pointsystems, respectively, where each $(\alpha,\beta)$ is generated by  $2{\rm {\it {\it rand}}}(1,2)-0.5$.

\begin{figure}[htbp]
\centerline{\includegraphics[height=8cm,width=14cm]{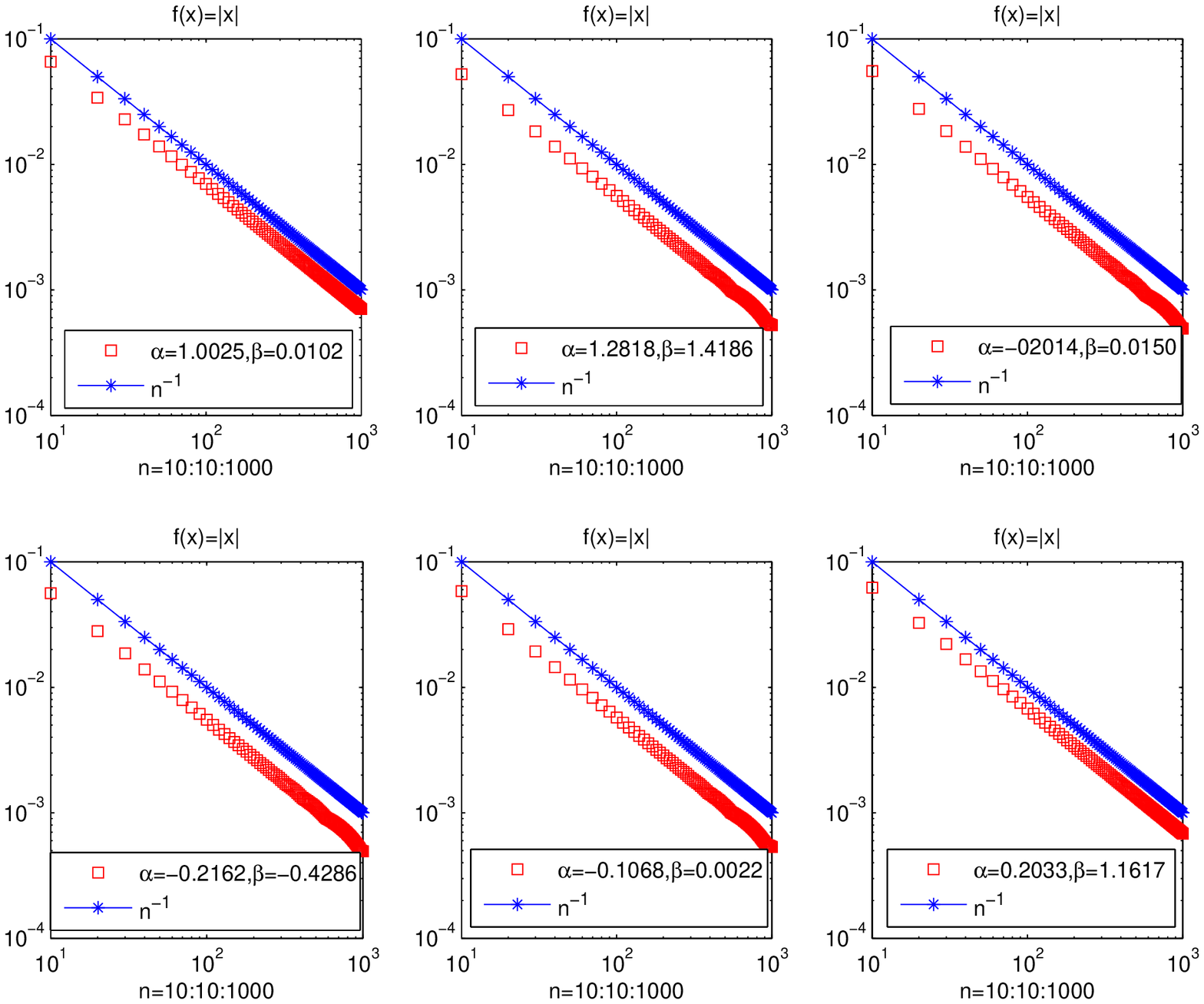}}
   \caption{$\max_{x=-1:0.001:1}|f(x)-L_n[f](x)|$ with $n=10:10:1000$ at the Jacobi-Gauss-Lobatto pointsystems for $f(x)=|x|$, respectively.}
\end{figure}

\begin{figure}[htbp]
\centerline{\includegraphics[height=8cm,width=14cm]{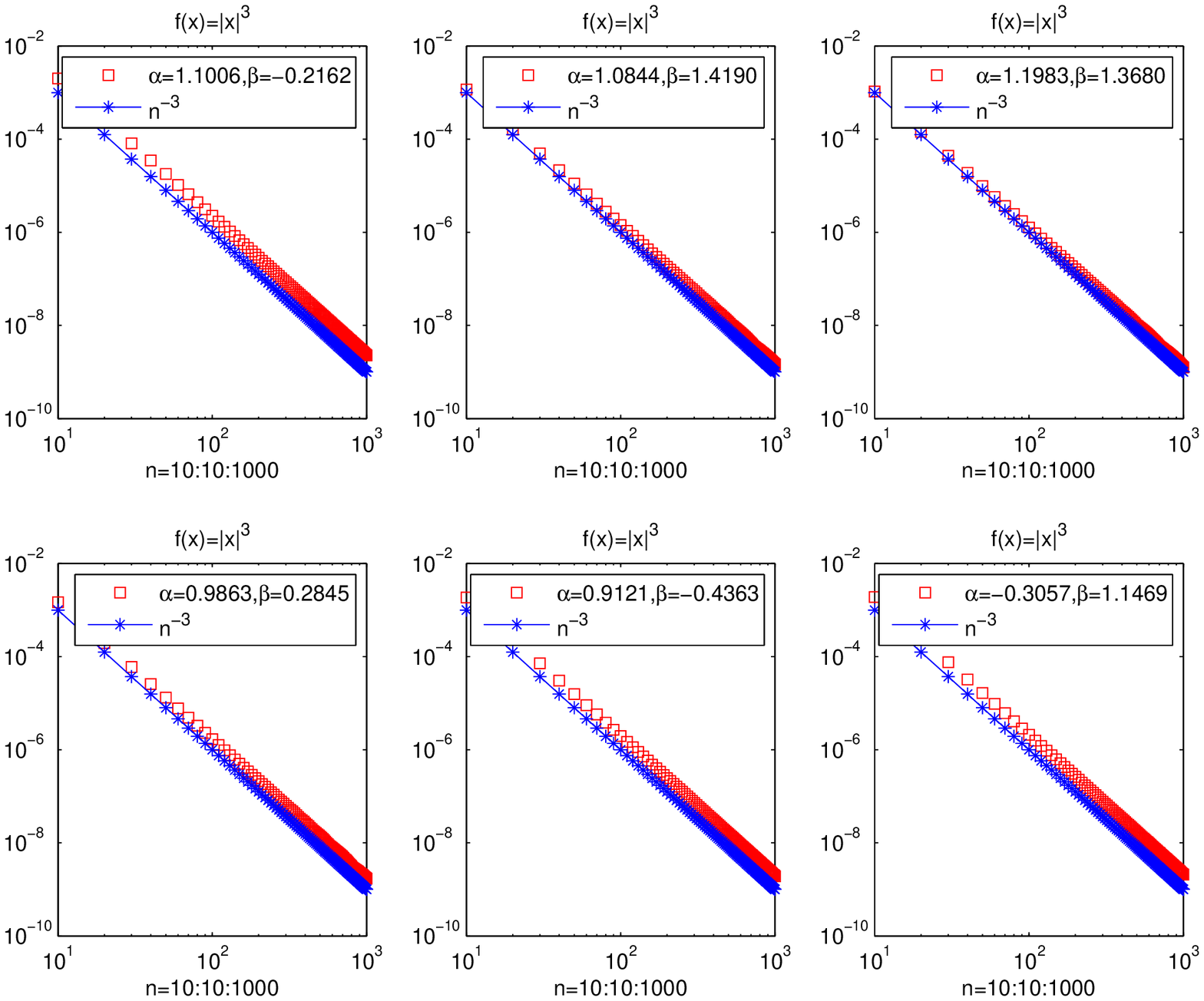}}
   \caption{$\max_{x=-1:0.001:1}|f(x)-L_n[f](x)|$ with $n=10:10:1000$ at the Jacobi-Gauss-Lobatto pointsystems for $f(x)=|x|^3$, respectively.}
\end{figure}

\subsection{General Jacobi-Gauss-Radau pointsystems}
Let
\begin{equation}\label{RadauGaussJacobi}
-1<x_n<x_{n-1}<\cdots<x_2<x_1<x_0=1
\end{equation}
be the roots of $(1-t)P_{n}^{(\alpha,\beta)}(t)=0$ ($\alpha,\beta>-1$), $x_k=\cos(\theta_k)$ and
$$
\omega(t)=(t-x_0)(t-x_1)\cdots(t-x_n), \quad \ell_k(t)=\frac{\omega(t)}{(t-x_k)\omega'(x_k)}.$$
Then
\begin{equation}
\ell_0(t)=\frac{P_{n}^{(\alpha,\beta)}(t)}{P_{n}^{(\alpha,\beta)}(1)},\quad\quad
{\displaystyle\ell_k(t)=\frac{(1-t)P_{n}^{(\alpha,\beta)}(t)}{(t-x_k)(1-x_k){P_{n}^{(\alpha,\beta)}}'(x_k)}},\quad k=1,2,\ldots,n.
\end{equation}
Additionally, for $t=\cos\theta\in [0,1]$, we have
$$
\ell_0(t)=\left\{\begin{array}{ll}O\left(\frac{n^{\alpha}}{n^{\alpha}}\right),&0\le \theta\le cn^{-1}\\
O\left(\theta^{-\alpha-\frac{1}{2}}n^{-\alpha-\frac{1}{2}}\right),&cn^{-1}\le\theta<\frac{\pi}{2}\end{array}\right.=O\left(n^{-\min\{0,\alpha+\frac{1}{2}\}}\right),
$$
and
\begin{equation}\begin{array}{lll}
\ell_k(t)&=&{\displaystyle -\frac{2\sin^2(\theta/2)P_{n}^{(\alpha,\beta)}(\cos\theta)}
{2\sin^2(\theta_k/2){P_{n}^{(\alpha,\beta)}}'(\cos\theta_k)2\sin((\theta-\theta_k)/2)\sin((\theta+\theta_k)/2)}}\\
&=&\left\{\begin{array}{ll}
O\left(n^{-\min\{0, \alpha+\frac{1}{2}\}}\right),&-1<\alpha\le -\frac{1}{2}\\
O\left(n^{-\min\{0,\frac{5}{2}-\alpha\}}\right),&\alpha>-\frac{1}{2}\end{array}\right.\end{array}
\end{equation}
following (4.27) similarly.

Similarly, for $t\in [-1,0]$, by $P_{n}^{(\beta,\alpha)}(-t)=(-1)^nP_{n}^{(\alpha,\beta)}(t)$, setting $t=-\cos\theta$ and $y_k=-x_{n-k+1}=\cos{\overline{\theta}_k}$ for $k=1,2,\ldots,n$, we obtain
$\ell_0(t)=O\left(n^{-\min\{\alpha+\frac{1}{2},\alpha-\beta\}}\right)$ and
\begin{equation}
\ell_k(t)=\left\{\begin{array}{ll}O\left(n^{-\min\{0,\alpha+\frac{1}{2},\alpha-\beta\}}\right),&-1<\alpha<\frac{1}{2}\\
O\left(n^{-\min\{0,\frac{1}{2}-\beta\}}\right),&\alpha\ge \frac{1}{2}\end{array}.\right.
\end{equation}
Thus for $t\in [-1,1]$, we get
\begin{equation}
\|\ell_k\|_{\infty}=\left\{\begin{array}{ll}O\left(n^{-\min\{0,\alpha+\frac{1}{2},\alpha-\beta\}}\right),&-1<\alpha\le\frac{1}{2}\\
O\left(n^{-\min\{0,\frac{1}{2}-\beta,\frac{5}{2}-\alpha,\alpha-\beta\}}\right),&\alpha>\frac{1}{2}\end{array}\right.
\end{equation}
for $k=0,1,2,\ldots,n$.

\begin{theorem}
Suppose $f(t)$ satisfies (2.9) and $\{x_j\}_{j=0}^{n}$  are the roots of $(1-t)P_n^{(\alpha,\beta)}(t)$ , then for $n\ge r+1$
\begin{equation}
E_n[f]=n^{-r}\cdot\left\{\begin{array}{ll}O\left(n^{-\min\{0,\alpha+\frac{1}{2},\alpha-\beta\}}\right),&-1<\alpha\le\frac{1}{2}\\
O\left(n^{-\min\{0,\frac{1}{2}-\beta,\frac{5}{2}-\alpha,\alpha-\beta\}}\right),&\alpha>\frac{1}{2}\end{array}\right.
\end{equation}
Particularly, we have for $(\alpha,\beta)\in \overline{S}$
\begin{equation}
\|E_n[f]\|_{\infty}=O\left(\frac{1}{n^{r}}\right),
\end{equation}
where $\overline{S}:= \left[-\frac{1}{2},\frac{5}{2}\right]\times \left(-1,\frac{1}{2}\right]- \left\{(\alpha,\beta):-\frac{1}{2}\le \alpha\le\frac{1}{2},\, \alpha<\beta\le \frac{1}{2}\right\}$.
\end{theorem}

Similarly we have

\begin{theorem}
Suppose $f(t)$ satisfies (2.9) and $\{x_j\}_{j=0}^{n}$  are the roots of $(1+t)P_n^{(\alpha,\beta)}(t)$ , then for $n\ge r+1$
\begin{equation}
E_n[f]=n^{-r}\cdot\left\{\begin{array}{ll}O\left(n^{-\min\{0,\beta+\frac{1}{2},\beta-\alpha\}}\right),&-1<\beta\le\frac{1}{2}\\
O\left(n^{-\min\{0,\frac{1}{2}-\alpha,\frac{5}{2}-\beta,\beta-\alpha\}}\right),&\beta>\frac{1}{2}\end{array}\right.
\end{equation}
Particularly, we have for $(\alpha,\beta)\in \widehat{S}$
\begin{equation}
\|E_n[f]\|_{\infty}=O\left(\frac{1}{n^{r}}\right),
\end{equation}
where $\widehat{S}:= \left(-1,\frac{1}{2}\right]\times \left[-\frac{1}{2},\frac{5}{2}\right]-\left\{(\alpha,\beta):-\frac{1}{2}\le \beta\le\frac{1}{2},\, \beta <\alpha\le\frac{1}{2}\right\}$.
\end{theorem}

{\sc Fig}. 4.11 shows the convergence rates for $f(x)=|x|$ at the Jacobi-Gauss-Radau pointsystems, where each $(\alpha,\beta)\in \overline{S}$ or $(\alpha,\beta)\in \widehat{S}$.

\begin{figure}[htbp]
\centerline{\includegraphics[height=5cm,width=6cm]{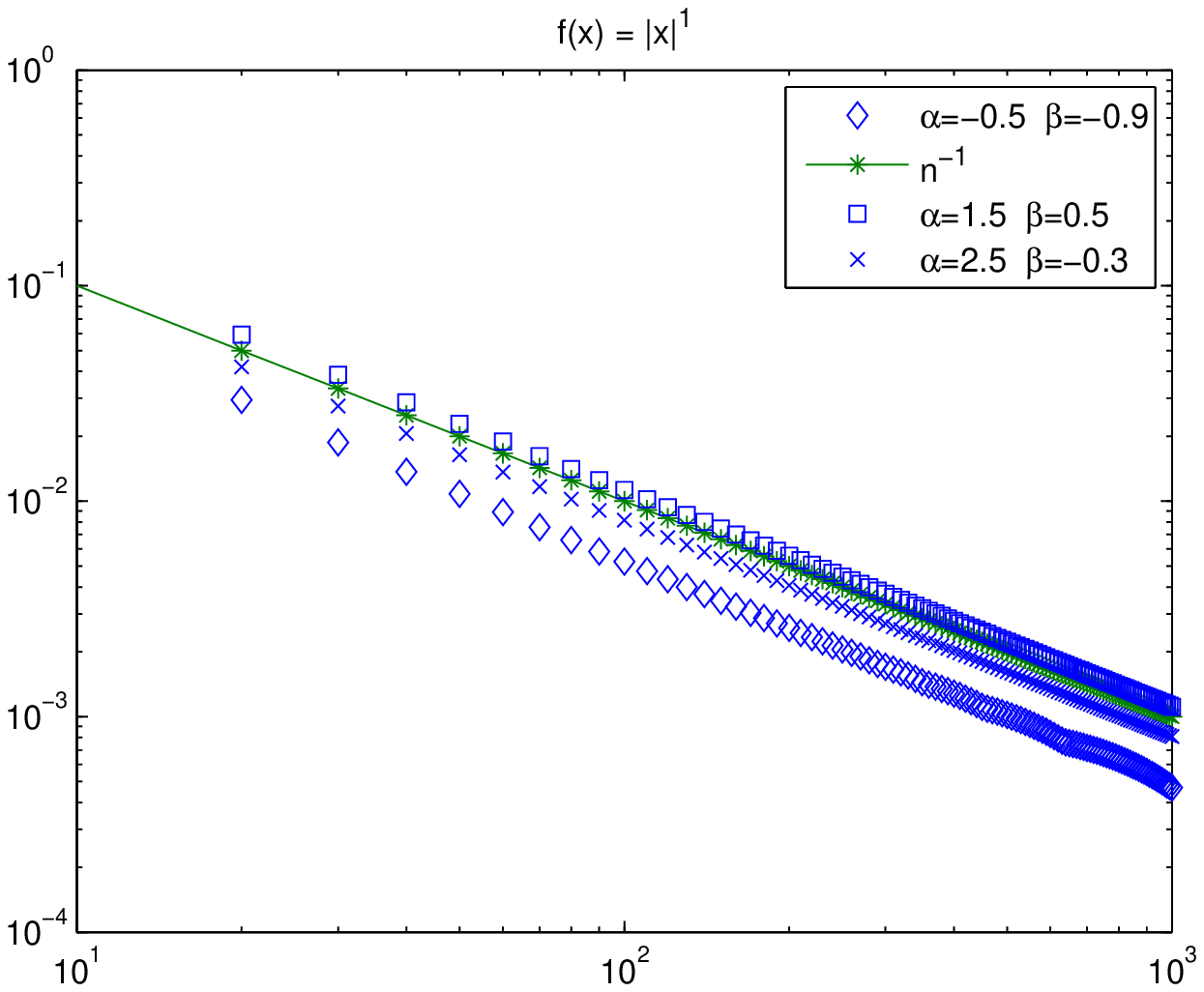}\hspace{1cm}\includegraphics[height=5cm,width=6cm]{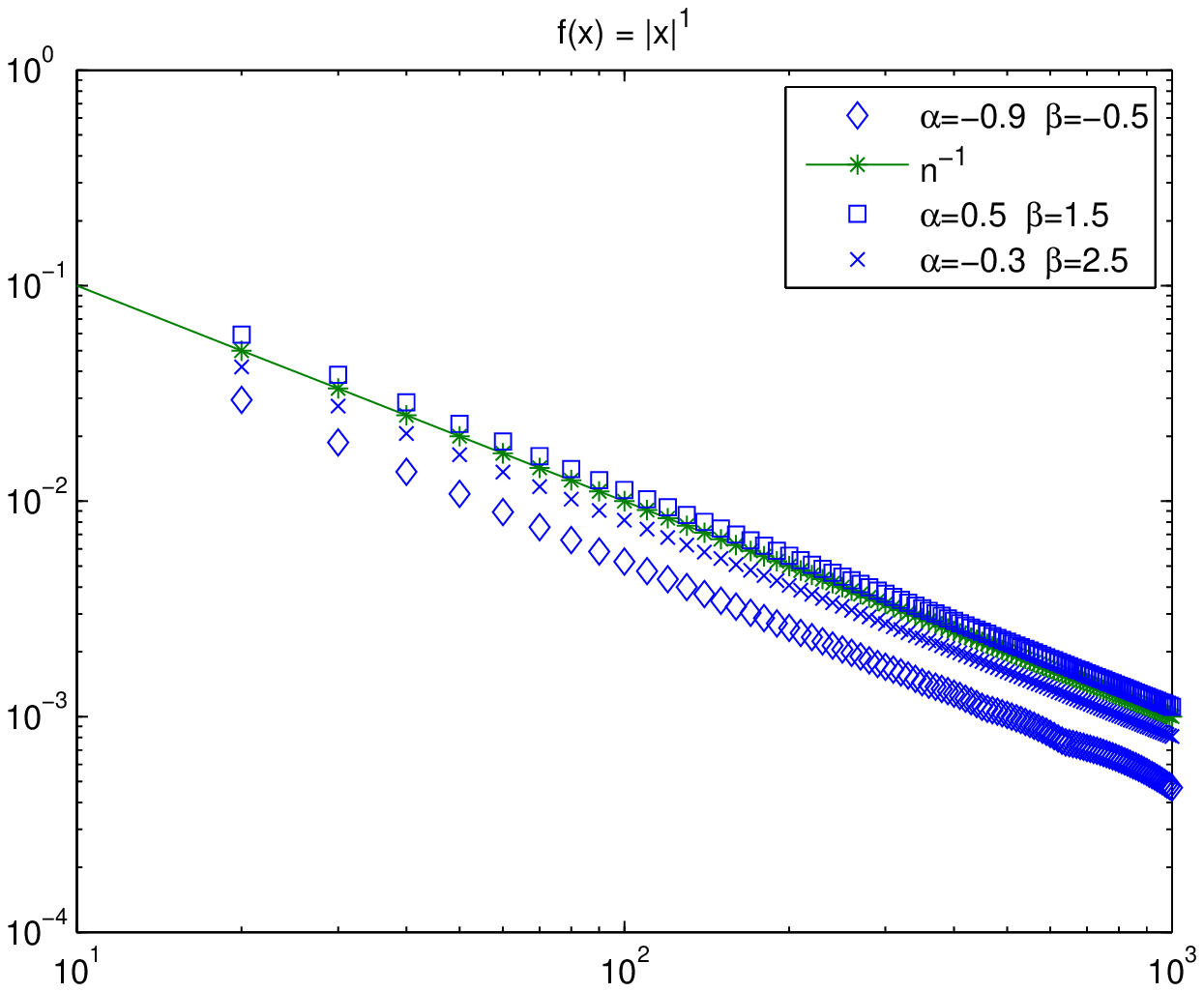}}
\caption{The absolute errors of  $\max_{x=-1:0.001:1}|f(x)-L_n[f](x)|$ for  Jacobi-Gauss-Radau pointsystems: the roots of $(1-t)P_n^{(\alpha,\beta)}(t)$ (left) and the roots of $(1+t)P_n^{(\alpha,\beta)}(t)$ (right) for $f(x)=|x|$, respectively.}
\end{figure}

\section{Final remarks}
The results in section 4 indicate the fact that the interpolations, for functions of limited regularities, at strongly normal pointsystems, Gauss-Jacobi pointsystems with $-1<\alpha,\beta\le\frac{1}{2}$, Jacobi-Gauss-Lobatto pointsystems with $(\alpha,\beta)\in S$, Gauss-Jacob-Radau pointsystems $\overline{S}$ or $\widehat{S}$, have the same convergence order  compared with the best polynomial approximation of the same degree. Numerical experiments give in line with the estimates.

In addition, numerical experiments also show that the same occurs for analytic or smooth functions. Here we illustrate the phenomenons by entire function $f(x)=e^x$, i.e., analytic throughout the complex plane, $f(x)=1/(1+25x^2)$, which is analytic in a neighborhood of $[-1, 1]$ but not throughout the complex plane, and $f(x)=e^{-1/x^2}$, which is not analytic in a neighborhood of $[-1, 1]$ but is infinitely differentiable in $[-1,1]$.

In {\sc Figs.} 5.1-5.3, the left columns are computed by zeros of Gauss-Jacobi polynomial $P^{(\alpha,\beta)}_n(x)$, the middles by Jacobi-Gauss-Lobatto $(1-x^2)P^{(\alpha,\beta)}_{n-2}(x)$, while the rights by Jacobi-Gauss-Radau $(1-x)P^{(\alpha,\beta)}_{n-1}(x)$ (first three cases) or  $(1+x)P^{(\alpha,\beta)}_{n-1}(x)$ (last three cases), respectively. From these figures, we see that the interpolations at these pointsystems including the Gauss-Legendre and Legendre-Gauss-Lobatto, achieve essentially the same approximation accuracy compared with those at the two Chebyshev piontsystems too.

\begin{figure}[htbp]
\centerline{\includegraphics[height=6.25cm,width=16cm]{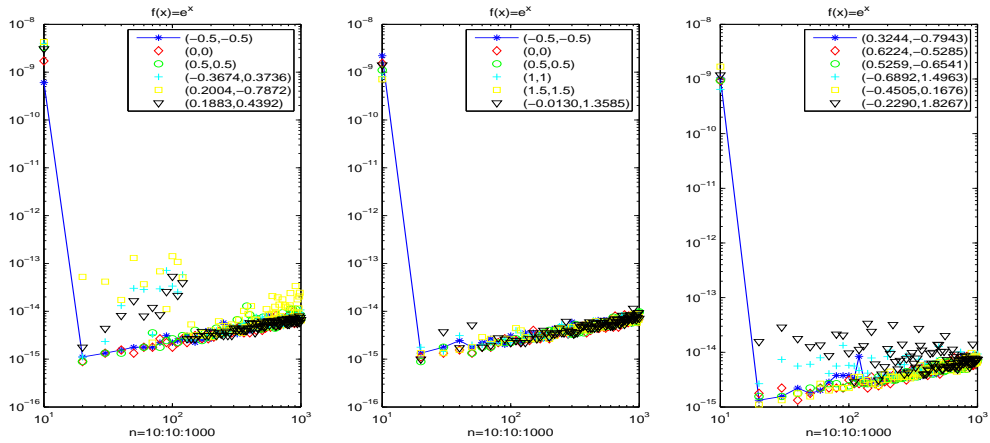}}
   \caption{$\max_{x=-1:0.001:1}|f(x)-L_n[f](x)|$ with $n=10:10:1000$ at Gauss-Jacobi pointsystems for $f(x)=e^x$, respectively.}
\end{figure}

\begin{figure}[htbp]
\centerline{\includegraphics[height=6.25cm,width=16cm]{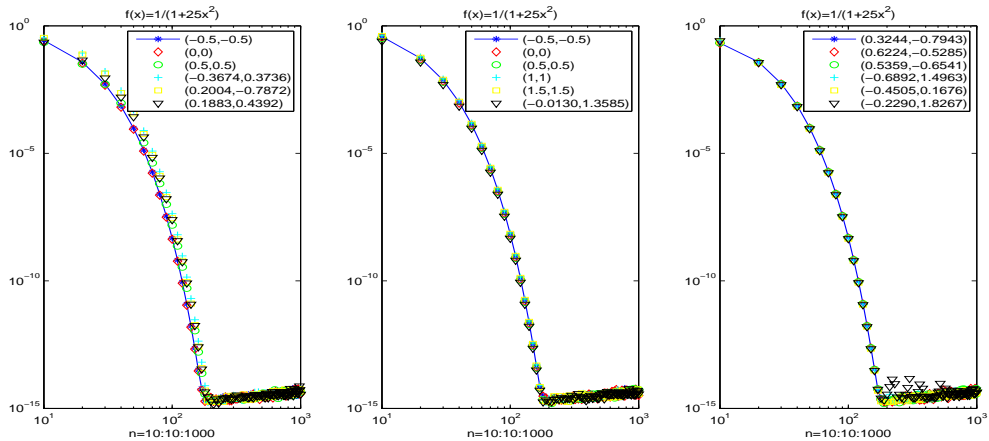}}
   \caption{$\max_{x=-1:0.001:1}|f(x)-L_n[f](x)|$ with $n=10:10:1000$ at Gauss-Jacobi pointsystems for $f(x)=\frac{1}{1+25x^2}$, respectively.}
\end{figure}

\begin{figure}[htbp]
\centerline{\includegraphics[height=6.25cm,width=16cm]{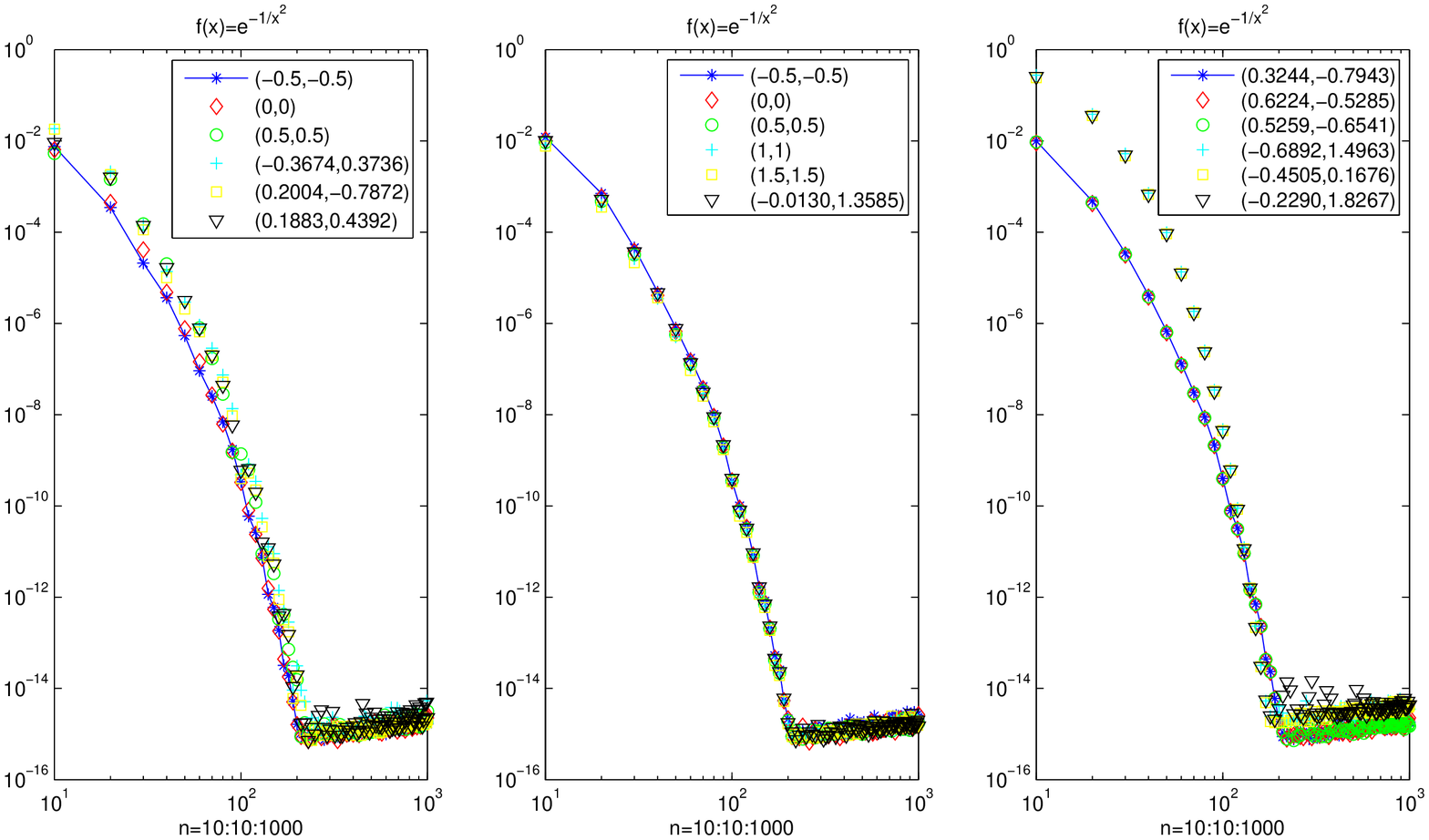}}
   \caption{$\max_{x=-1:0.001:1}|f(x)-L_n[f](x)|$ with $n=10:10:1000$ at Jacobi-Gauss-Radau pointsystems for $f(x)=e^{-1/x^2}$, respectively.}
\end{figure}

It is interesting to noting that the interpolation approximation polynomial $L_n[f]$ challenges the best approximation polynomial $p^*_{n-1}$ of $f$ at the above nice pointsystems not only on the equally asymptotic order on the convergence rate, but also on the faster convergence on the first derivative or second derivative approximation by the polynomials, which has plenty applications in spectral methods \cite{Hesthaven,Trefethenbook}.

{\sc Fig.} 5.4 illustrates that the convergence rates $\max_{x=-1:0.001:1}|f'(x)-[p^*_{n-1}]'(x)|$ and  $\max_{x=-1:0.001:1}|f''(x)-[p^*_{n-1}]''(x)|$ reduces 2-order and 3-order, respectively, compared with $\max_{x=-1:0.001:1}|f(x)-p_{n-1}^*(x)|$ for $f(x)=|x|^5$ or $|x|^7$. Here, the best approximation polynomial $p^*_{n-1}$ is obtained by {\bf remez} algorithm in {\sc Chebfun} system \cite{Chebfun}.

However, the interpolation polynomial $L_n[f]$ at the above nice pointsystems performs much better than $p^*_{n-1}$ for approximation $f'$ and $f''$ by $L_n'[f]$, $L_n''[f]$, $[p^*_{n-1}]'$ and $[p^*_{n-1}]''$, respectively. Here, we use usual pointsystems to show the performs (see {\sc Figs.} 5.5-5.7).

\begin{figure}[htbp]
\centerline{\includegraphics[height=5cm,width=7cm]{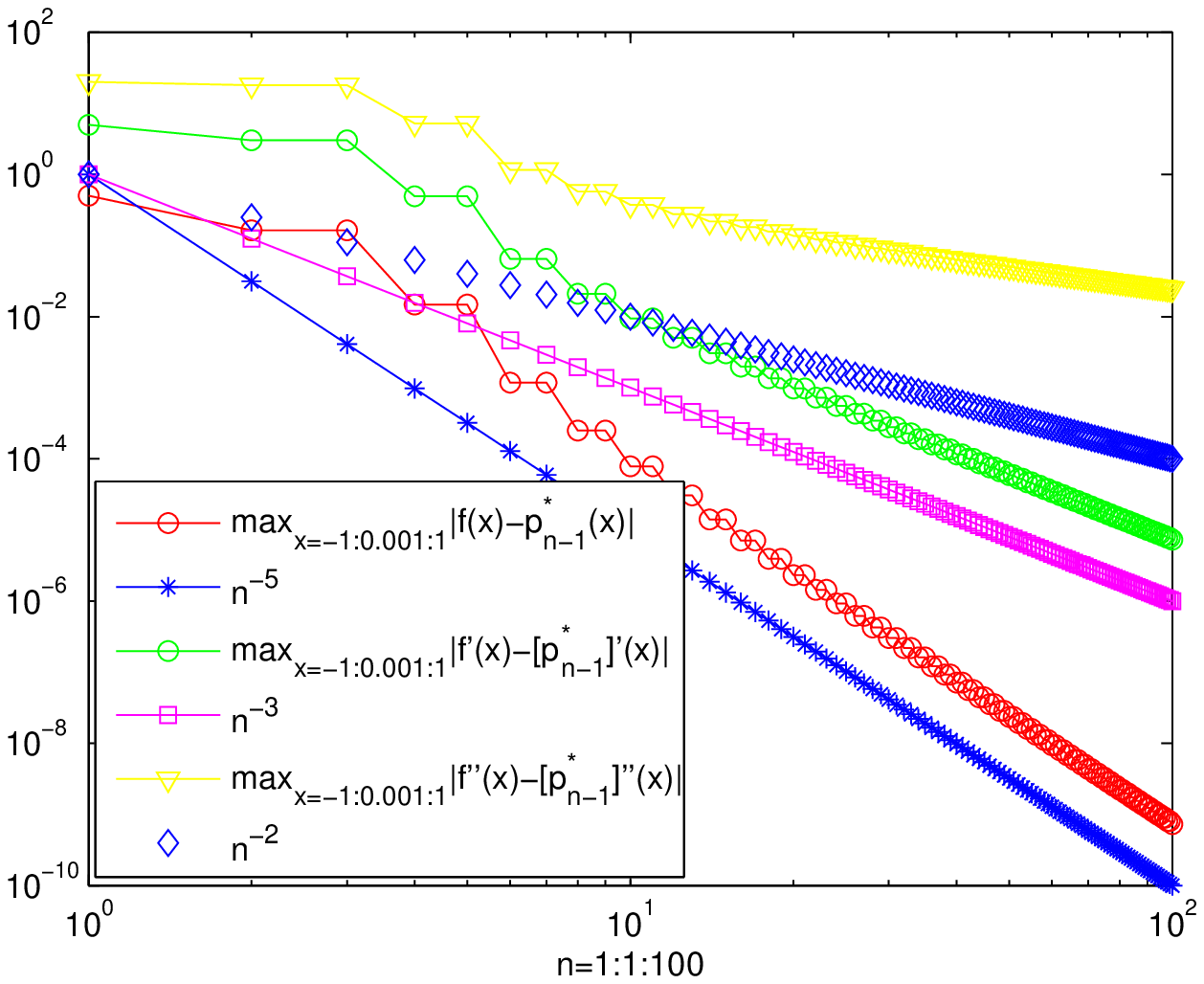}\includegraphics[height=5cm,width=7cm]{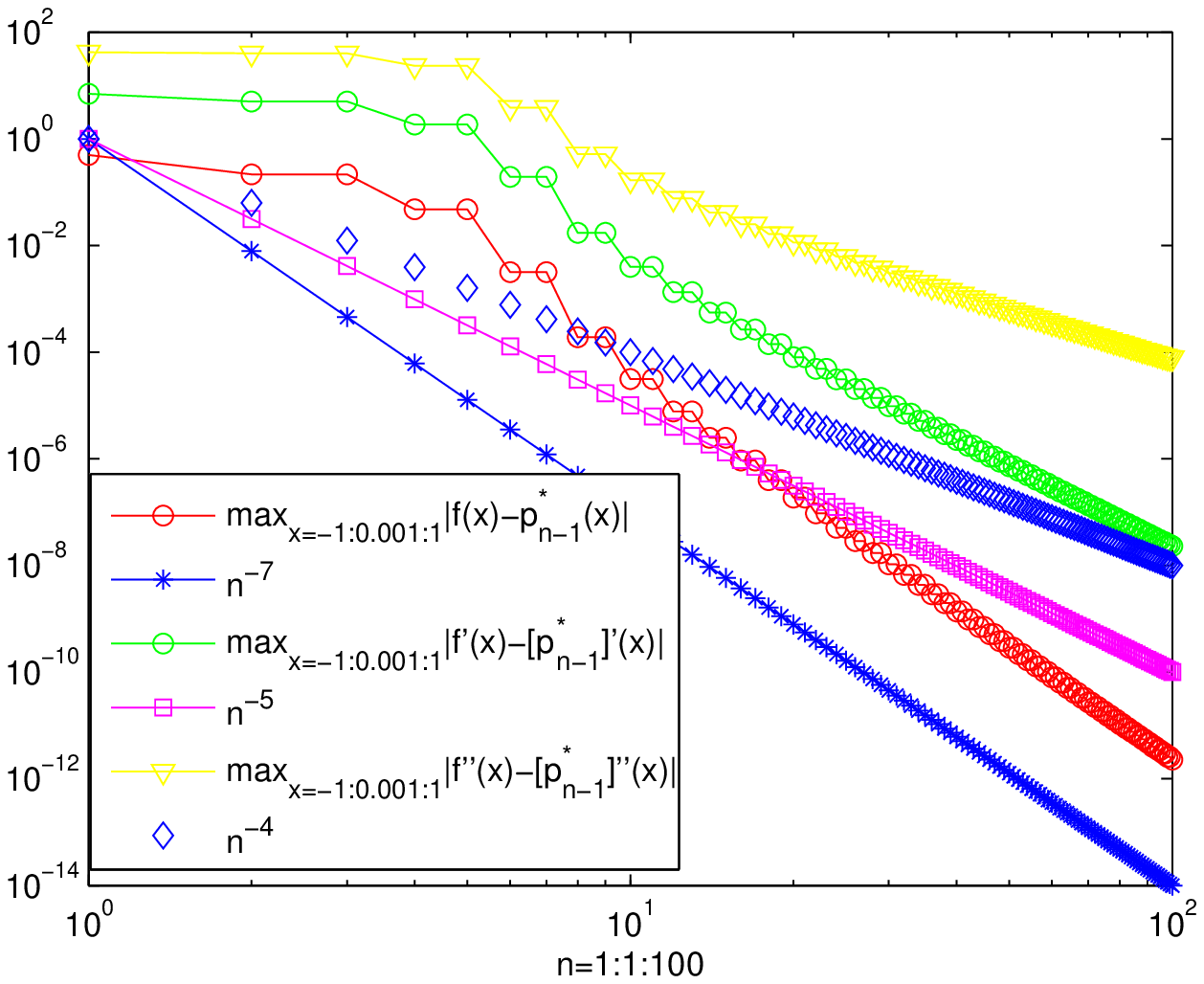}}
   \caption{$\max_{x=-1:0.001:1}|f^{(m)}(x)-[p^*_{n-1}]^{(m)}(x)|$ with $n=10:1:100$  for $f(x)=|x|^5$ (left), $f(x)=|x|^7$ (right) and $m=0,1,2$, respectively.}
\end{figure}

\begin{figure}[htbp]
\centerline{\includegraphics[height=6.25cm,width=16cm]{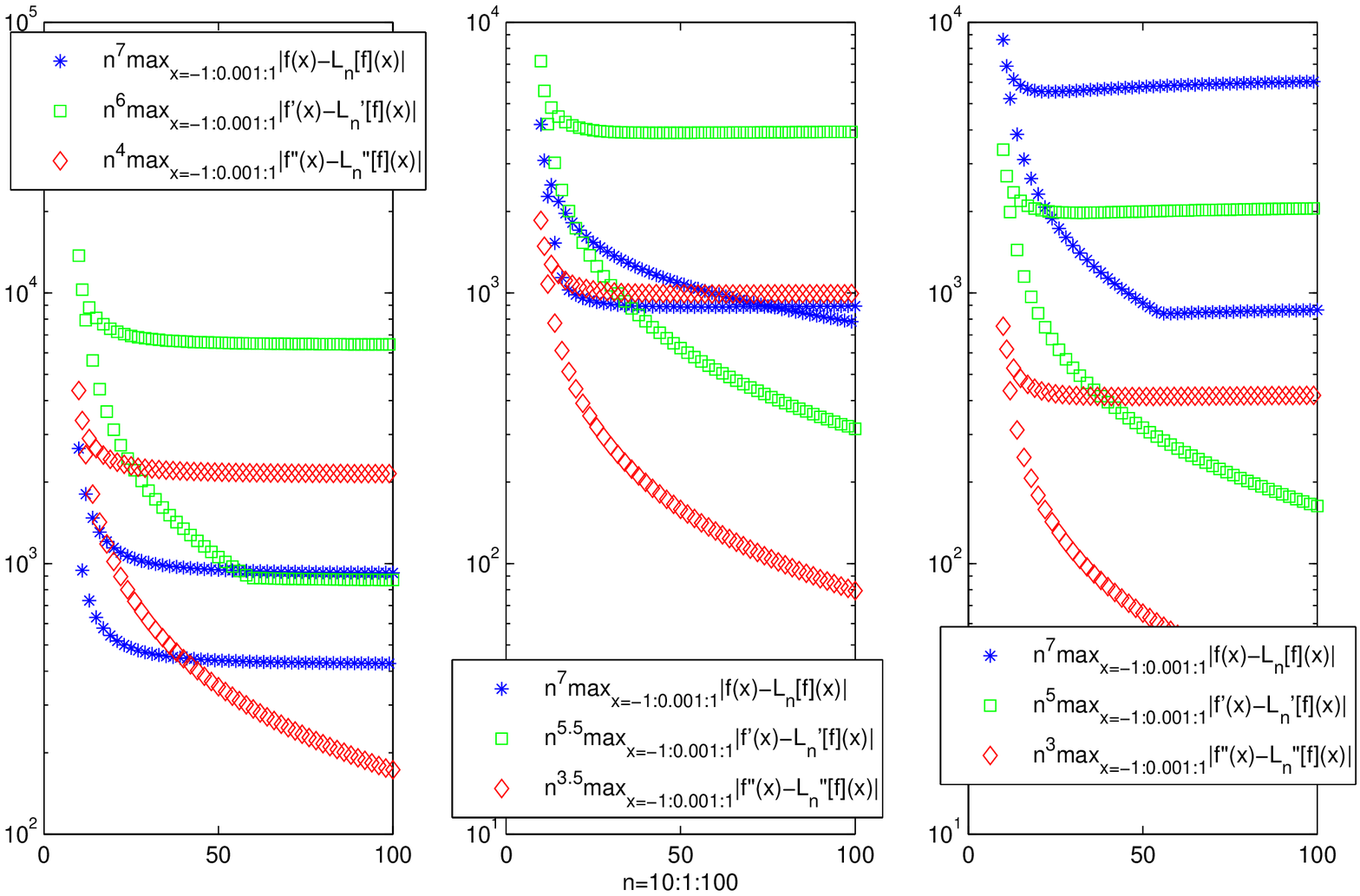}}
   \caption{$n^{s(m)}\max_{x=-1:0.001:1}|f^{(m)}(x)-L_n^{(m)}[f](x)|$ with $n=10:1:100$ at Gauss-Jacobi pointsystems for $f(x)=|x|^7$ and $m=0,1,2$,  and (-0.5,-0.5) (left), (0,0) (middle), (0.5,0.5) (right), respectively. }
\end{figure}

\begin{figure}[htbp]
\centerline{\includegraphics[height=6.25cm,width=16cm]{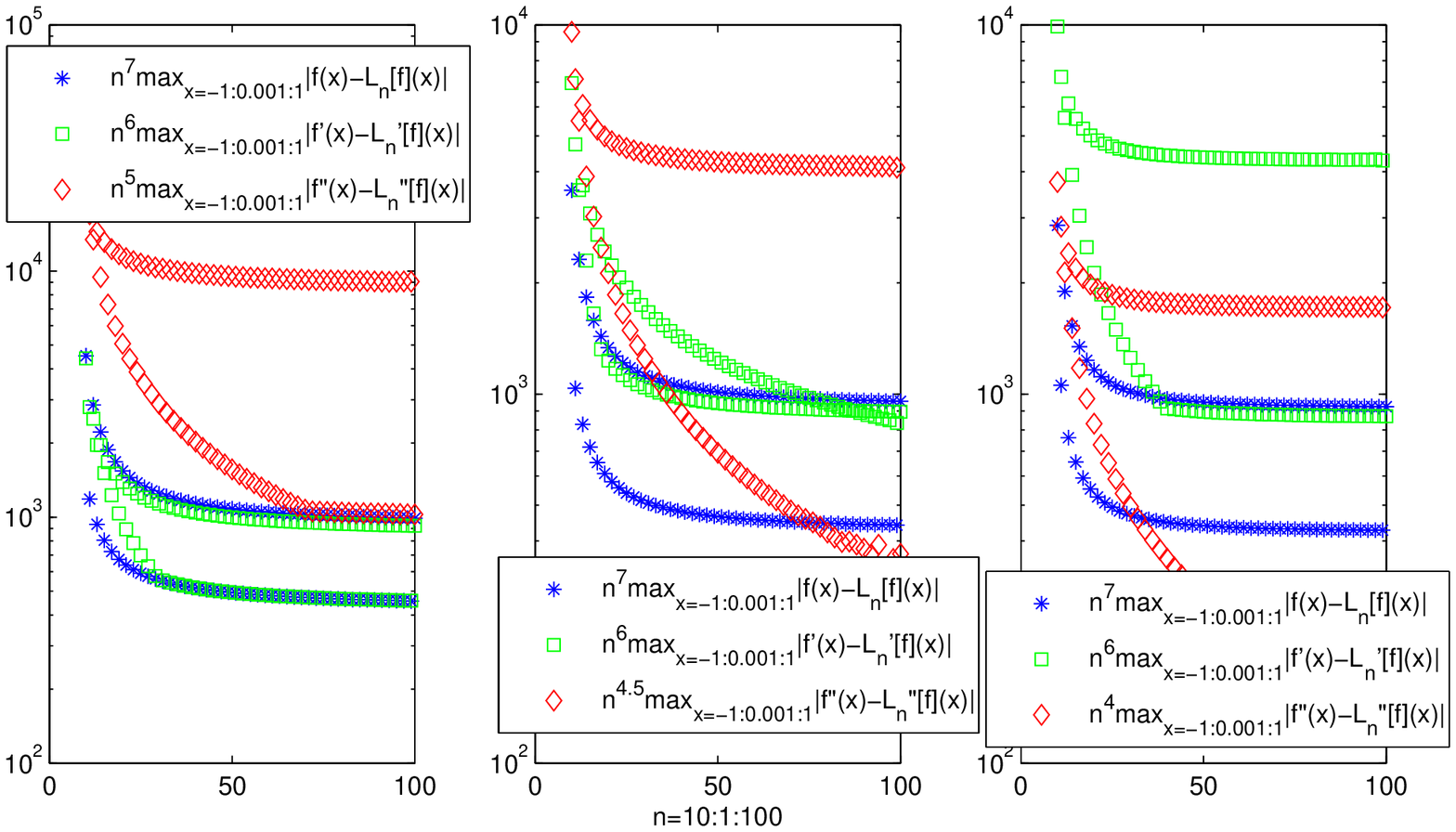}}
   \caption{$n^{s(m)}\max_{x=-1:0.001:1}|f^{(m)}(x)-L_n^{(m)}[f](x)|$ with $n=10:1:100$ at Jacobi-Gauss-Lobatto pointsystems for $f(x)=|x|^7$ and $m=0,1,2$, and (0.5,0.5) (left), (1,1) (middle), (1.5,1.5) (right), respectively.}
\end{figure}

\begin{figure}[htbp]
\centerline{\includegraphics[height=6.25cm,width=16cm]{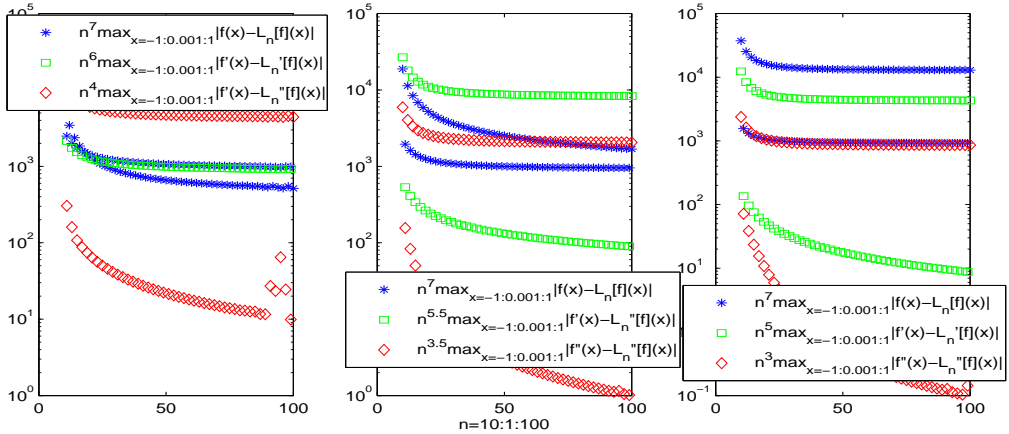}}
   \caption{$n^{s(m)}\max_{x=-1:0.001:1}|f^{(m)}(x)-L_n^{(m)}[f](x)|$
    with $n=10:1:100$ at Jacobi-Gauss-Radau pointsystems with $x_0=1$ for $f(x)=|x|^7$ and $m=0,1,2$, and (-0.5,-0.5) (left), (0,0) (middle), (0.5,0.5) (right),  respectively.}
\end{figure}

{\sc Fig.} 5.8 shows that the convergence rates $\max_{x=-1:0.001:1}|f^{(m)}(x)-[p^*_{n-1}]^{(m)}(x)|$ and  $\max_{x=-1:0.001:1}|f^{(m)}(x)-L_n^{(m)}[f](x)|$ for $f(x)=|x|^5$ and $|x|^7$, respectively, where $L_n[f]$ is the interpolation at Jacobi-Gauss-Lobatto points $\{x_k=\cos\left(\frac{k\pi}{n-1}\right)\}_{k=0}^{n-1}$.
\begin{figure}[htbp]
\centerline{\includegraphics[height=5cm,width=7cm]{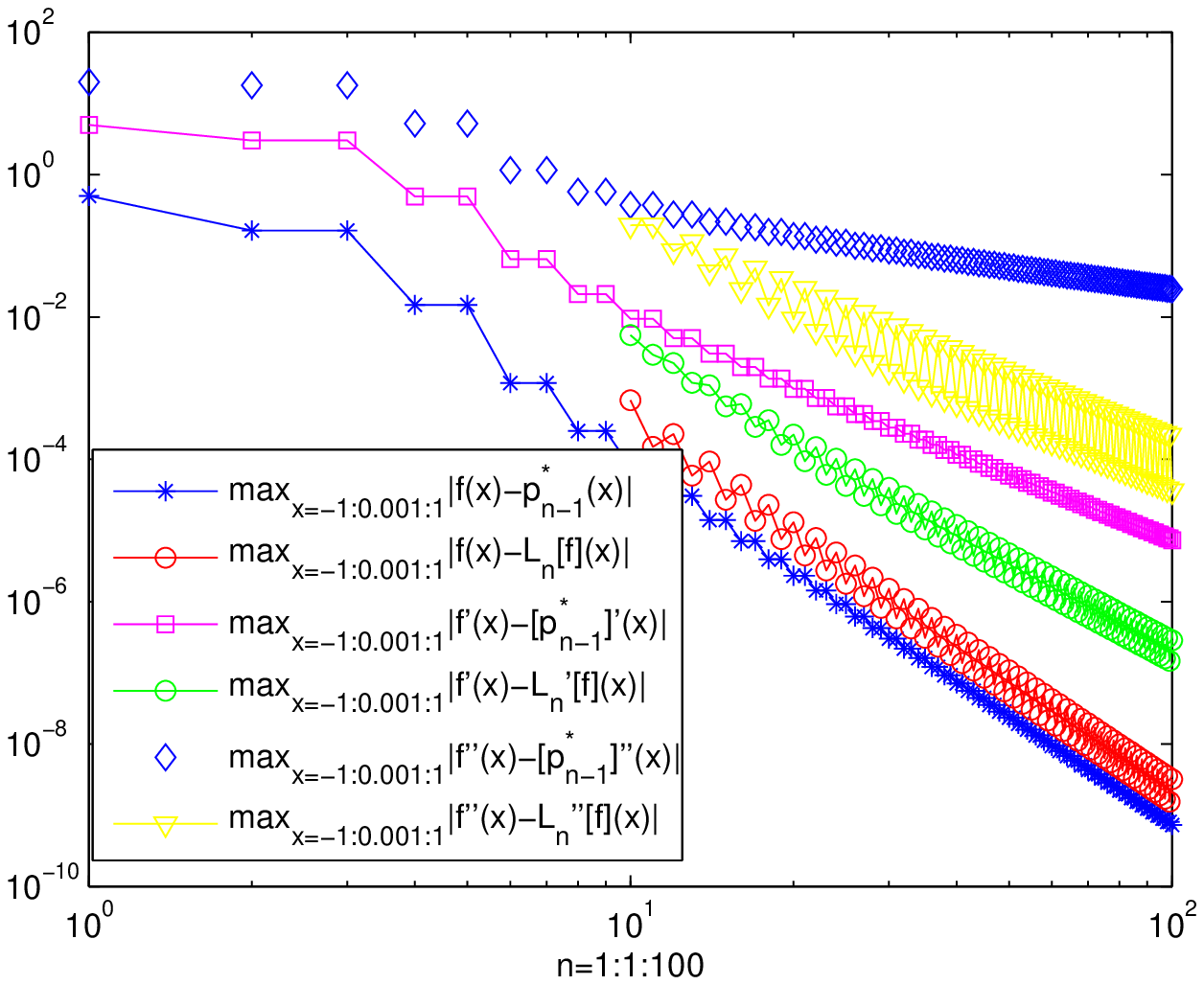}\includegraphics[height=5cm,width=7cm]{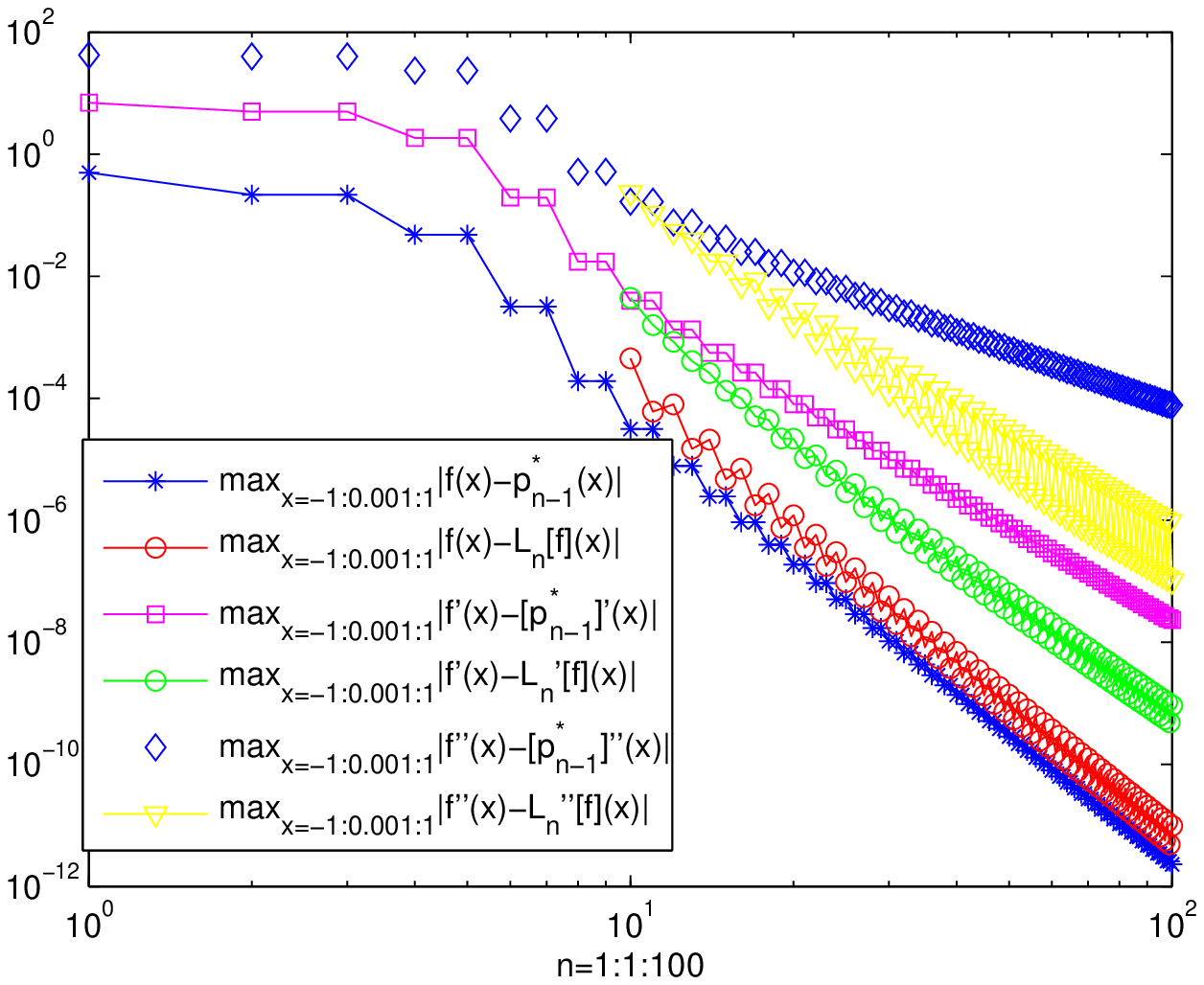}}
   \caption{The convergence rates $\max_{x=-1:0.001:1}|f^{(m)}(x)-[p^*_{n-1}]^{(m)}(x)|$ and  $\max_{x=-1:0.001:1}|f^{(m)}(x)-L_n^{(m)}[f](x)|$ for $f(x)=|x|^5$ (left) and $|x|^7$ (right) with $m=0,1,2$, respectively, where $L_n[f]$ is the interpolation at Jacobi-Gauss-Lobatto points $\{x_k=\cos\left(\frac{k\pi}{n-1}\right)\}_{k=0}^{n-1}$.}
\end{figure}

\vspace{0.1cm} {\bf Acknowledgement}. The author is grateful to
Professor Kelzon for his kind help, and to  Professor Yuri Wainerman
 for her sending me three pages of her Ph. D thesis including the interesting Lemma 3.1.  The author thanks  Chaoxu Pei at Florida State University
and Yulong Lu at University of Warwick   for their checking up every detail and constructive comments. The author also thanks Dr. Guo He and  Guidong Liu, at Central South University, for their helpful and insightful discussion on the convergence rates of the derivatives of interpolations and estimats on the Jacobi-Gauss-Radau pointsystems.

\end{document}